\documentclass[11pt]{article}
\usepackage[matrix,arrow,curve,cmtip]{xy}
\usepackage{amssymb}
\usepackage{latexsym}
\usepackage{theorem}
\usepackage[a4paper,width=140mm,height=220mm]{geometry}
\usepackage{titlesec}
\titleformat{\section}[hang]%
{\scshape\filcenter\large}{\oldstylenums{\thesection}.}{1ex}{}%
\renewcommand{\thesubsection}{\Alph{subsection}}
\titleformat{\subsection}[runin]%
{\bfseries\normalsize}{\thesubsection.}{1ex}{}[.\hspace{1ex}]

\def\to{\mbox{$\xymatrix@1@C=5mm{\ar@{->}[r]&}$}}
\def\tto{\mbox{$\xymatrix@1@C=5mm{\ar@{=>}[r]&}$}}
\def\halfcirc{\begin{picture}(0,0)\put(0,3){\oval(2,2)[l]}\end{picture}}
\def\incl{\mbox{$\xymatrix@1@C=5mm{\ar@{->}[r]|<{\halfcirc}&}$}}
\def\distsign{\begin{picture}(0,0)\put(0,0){\circle{4}}\end{picture}}
\def\dist{\mbox{$\xymatrix@1@C=5mm{\ar@{->}[r]|{\distsign}&}$}}
\def\biar{\mbox{$\xymatrix@1@C=5mm{\ar@<1.5mm>[r]\ar@<-0.5mm>[r]&}$}}
\def\bidist{\mbox{$\xymatrix@1@C=5mm{\ar@<1.5mm>[r]|{\distsign}\ar@<-0.5mm>[r]|{\distsign}&}$}}
\def\iso{\mbox{$\xymatrix@1@C=6mm{\ar@{->}[r]^{\sim}&}$}}

\def\inlineadj#1#2{\mbox{$\xymatrix@C=15mm{\ar@{}[r]|{\bot}\ar@<1mm>@/^2mm/[r]^{{#1}} & \ar@<1mm>@/^2mm/[l]^{{#2}}}$}}

\def\inlinerevadj#1#2{\mbox{$\xymatrix@C=15mm{\ar@{}[r]|{\bot}\ar@<-1mm>@/_2mm/[r]_{{#1}} & \ar@<-1mm>@/_2mm/[l]_{{#2}}}$}}

\def\inlineadjretr#1#2{\mbox{$\xymatrix@C=15mm{\ar@{}[r]|{\bot}\ar@{>->}@<-1mm>@/_2mm/[r]_{{#1}} & \ar@{->>}@<-1mm>@/_2mm/[l]_{{#2}}}$}}

\def\inlineadjdist#1#2{\mbox{$\xymatrix@C=15mm{\ar@{}[r]|{\bot}\ar@<1mm>@/^2mm/[r]|{\distsign}^{{#1}} & \ar@<1mm>@/^2mm/[l]|{\distsign}^{{#2}}}$}}

\newcommand\arr[2][7mm]{\mbox{$\xymatrix@1@C= #1 {\ar@{->}^-{{#2}}[r]&}$}}

\newcommand\twocell[2][7mm]{\mbox{$\xymatrix@1@C= #1 {\ar@{=>}^{{#2}}[r]&}$}}

\def\inlineeqvl#1#2{\mbox{$\xymatrix@C=15mm{\ar@<-.5mm>@{}[r]|{\sim}\ar@<1mm>@/^2mm/[r]^{{#1}} & \ar@<1mm>@/^2mm/[l]^{{#2}}}$}}

\def\upincl\mbox{$\xymatrix@R=10mm{ \\ \ar@{^{(}->}[u]}$}

\def\endoar#1#2{\mbox{\xymatrix@1{{#1}\ar@(u,r)|{#2}}}}

\def\down{\begin{picture}(6,0)
\put(0,0){$\downarrow$}
\end{picture}}

\def\endoar#1#2{\mbox{\xymatrix{{#1}\ar@(u,r)|{#2}}}}

\def\mono{\mbox{$\xymatrix@1@C=4.5mm{~\ar@{>->}[r]&}$}}

\newtheorem{theorem}{Theorem}[section]
\newtheorem{lemma}[theorem]{Lemma}
\newtheorem{definition}[theorem]{Definition} 
\newtheorem{proposition}[theorem]{Proposition}
\newtheorem{corollary}[theorem]{Corollary}
{\theorembodyfont{\upshape}\newtheorem{example}[theorem]{Example}}
{\theorembodyfont{\upshape}\newtheorem{remark}[theorem]{Remark}}
\newcommand{\proof}{\noindent {\em Proof\ }: }
\def\endofproof{$\mbox{ }\hfill\Box$\par\vspace{1.8mm}\noindent}

\def\eqref#1{(\ref{#1})}
\def\End{{\sf End}}

\def\Rel{{\sf Rel}}
\def\Frm{{\sf Frm}}
\def\LH{{\sf LH}}
\def\et{_{\sf \acute{e}t}}

\def\slh{_{\sf slh}}

\def\o{^{\sf o}}
\def\Loc{{\sf Loc}}
\def\boldalpha{\mbox{\boldmath$\alpha$}}
\def\boldbeta{\mbox{\boldmath$\beta$}}
\def\bigmid{\ \Big|\ }

\def\lpg{_{\sf lpg}}
\def\lpr{_{\sf lpr}}
\def\pg{_{\sf pg}}
\def\pr{_{\sf pr}}
\def\ta{_{\sf ta}}
\def\c{_{\sf c}}
\def\Sh{{\sf Sh}}
\def\si{_{\sf si}}
\def\Mod{{\sf Mod}}
\def\etal{~{\it et~al.}}
\def\Cocont{{\sf Cocont}}
\def\Ord{{\sf Ord}}
\def\:{\colon}
\def\impl{\Rightarrow}
\def\2{{\bf 2}}

\def\Set{{\sf Set}}
\def\QUANT{{\sf QUANT}}

\def\cc{_{\sf cc}}
\def\op{^{\sf op}}

\def\Sup{{\sf Sup}}
\def\Cat{{\sf Cat}}
\def\Dist{{\sf Dist}}
\def\Map{{\sf Map}}
\def\id{{\sf id}}
\def\Q{{\cal Q}}

\def\P{{\cal P}}

\def\F{{\cal F}}
\def\G{{\cal G}}

\def\colim{\mathop{\rm colim}}
\def\lim{\mathop{\rm lim}}
\def\sup{{\sf sup}}
\def\bbA{\mathbb{A}}
\def\bbB{\mathbb{B}}
\def\bbC{\mathbb{C}}
\def\bbD{\mathbb{D}}

\def\tensor{\otimes}

\def\coop{^{\sf coop}}

\def\slice#1{%
\setlength{\unitlength}{1ex}
\begin{picture}(2.7,1)(0,0)
\put(-0.2,-0.2){\mbox{$/$}}
\put(0.5,-0.7){\mbox{${#1}$}}
\end{picture}
\setlength{\unitlength}{1pt}}

\def\raiseletter#1{
\setlength{\unitlength}{1ex}
\begin{picture}(1,1)(0,0)
\put(0,0.7){\mbox{${#1}$}}
\end{picture}
\setlength{\unitlength}{1pt}}

\title{On principally generated quantaloid-modules in general,\\
and skew local homeomorphisms in particular}
\author{Hans Heymans\footnote{Department of Mathematics and Computer Science, University of Antwerp, Middelheimlaan 1, 2020 Antwerpen, Belgium, {\tt hans.heymans@ua.ac.be}} \ and Isar Stubbe\footnote{Postdoctoral Fellow of the Research Foundation Flanders (FWO), Department of Mathematics and Computer Science, University of Antwerp, Middelheimlaan 1, 2020 Antwerpen, Belgium, {\tt isar.stubbe@ua.ac.be}}}
\date{Written: December 18, 2007 \\ Submitted: March 5, 2008 \\ Revised: March 31, 2009}

\begin{document}

\maketitle

\begin{abstract}
Ordered sheaves on a small quantaloid $\Q$ have been defined in terms of $\Q$-enriched categorical structures; they form a locally ordered category $\Ord(\Q)$. The free-cocompletion KZ-doctrine on $\Ord(\Q)$ has $\Mod(\Q)$, the quantaloid of $\Q$-modules, as category of Eilenberg-Moore algebras. In this paper we give an intrinsic description of the Kleisli algebras: we call them the {\em locally principally generated $\Q$-modules}. We deduce that $\Ord(\Q)$ is biequivalent to the 2-category of locally principally generated $\Q$-modules and left adjoint module morphisms. The example of locally principally generated modules on a locale $X$ is worked out in full detail: relating $X$-modules to objects of the slice category $\Loc\slice{X}$, we show that ordered sheaves on $X$ correspond with {\em skew local homeomorphisms into $X$} (like sheaves on $X$ correspond with local homeomorphisms into $X$).
\end{abstract}

\section{Introduction}\label{A}

\subsection*{Locales and quantales, sheaves and logic}
A {\em locale} $X$ is a complete lattice in which finite infima distribute over arbitrary suprema. A particular class of examples of locales comes from topology: the open subsets of any topological space form a locale. But not every locale arises in this way, whence the slogan that locales are ``pointfree topologies'' [Johnstone, 1983]. There is a ``pointfree'' way to do sheaf theory: a sheaf $F$ on a locale $X$ is a functor $F\:X\op\to\Set$ satisfying gluing conditions. The collection of all such functors, together with natural transformations between them, forms the topos $\Sh(X)$ of {\em sheaves on $X$}. One of the many close ties between logic and sheaf theory, which is of particular interest to us, is that the {\em internal logic} of $\Sh(X)$ is an intuitionistic higher-order predicate logic with $X$ as object of truth values [Mac Lane and Moerdijk, 1991; Borceux, 1994; Johnstone, 2002]. To borrow a phrase from [Reyes, 1977] and others, sheaf theory thus serves as {\em algebraic logic}.
\par
The definition of locale can be restated: $X$ is a complete lattice and $(X,\wedge,\top)$ is a monoid 
such that the multiplication distributes on both sides over arbitrary suprema. It is natural to generalise this: a {\em quantale} $Q=(Q,\circ,1)$ is, by definition, a monoid structure on a complete lattice such that the multiplication distributes on both sides over arbitrary suprema [Mulvey, 1986; Rosenthal, 1992]. Because the monoid structure of a locale is obviously commutative, but the one for a quantale need not be, one can think of quantales as ``pointfree non-commutative topologies''. Examples of quantales, other than locales, can be found in algebra and geometry [Mulvey and Pelletier, 2001; Resende, 2007], in logic [Yetter, 1990], in computer science [Abramsky and Vickers, 1993; Rosicky, 200]. In the spirit of (enriched) category theory [Kelly, 1982], it is not hard to see that a quantale is precisely a monoid in the symmetric monoidal closed category $\Sup$ of complete lattices and suprema-preserving functions. And a {\em quantaloid} $\Q$ is then defined as a category enriched in $\Sup$ (so a quantale is the same thing as a quantaloid with one object, precisely as a group is a groupoid with one object) [Rosenthal, 1996].
\par
The success of sheaf theory to study logic from locales, and the useful generalisation from locales to quantales (and even quantaloids), make one wonder about the ``logic of sheaves on quantales''. However, it is not at all straightforward to define ``sheaves on a quantale''! Many different definitions have been proposed by many different authors, e.g.\ [Borceux and Van den Bossche, 1986; Mulvey and Nawaz, 1995; Gylys 2001; Garraway, 2005], but often only for particular classes of quantales. In previous work we have taken the following stance on the matter: whereas sheaves on a locale $X$ can be described in terms of sets equipped with an $X$-valued equality relation [Lawvere, 1973; Fourman and Scott, 1979; Borceux, 1994], the non-commutativity of the quantale multiplication forces sheaves on a quantale $Q$ to be sets equipped with a $Q$-valued {\em inequality} relation\footnote{We make a remark about $\Q$-valued {\em equalities} at the end of this Introduction.} [Borceux and Cruciani, 1998; Stubbe, 2005b]. In other words, our attention goes to the category of {\em ordered sheaves on a quantale} (or even quantaloid), which we see as ``algebraic non-commutative logic''.
\par
More precisely, Stubbe [2005b] studied ordered sheaves on a quantaloid $\Q$ in terms of $\Q$-enriched categories, thus generalising to the non-commutative case the work of [Walters, 1981; Borceux and Cruciani, 1998] on locales. In this paper we shall show that they can equivalently be described as particular $\Q$-modules. (If $\Q$ is a quantaloid, then a $\Q$-module is by definition a $\Sup$-enriched functor $F\:\Q\op\to\Sup$. For a quantale $Q$ this reduces to a complete lattice on which $Q$ acts.) This in turn can be applied to a locale $X$, and we find a characterisation of the relevant $X$-modules as particular locale morphisms with codomain $X$. We speak of {\em principally generated $\Q$-modules} in general, and {\em skew local homeomorphisms into $X$} in particular, as we shall introduce next.

\subsection*{Principally generated $\Q$-modules}

To introduce this novel notion in $\Q$-module theory, around which this article is centered, we first recall a simple fact from order theory: The well-known adjunction between the category $\Ord$ of ordered sets and order-preserving functions on the one hand, and the category $\Sup$ of complete lattices and supremum-preserving functions on the other,
\begin{equation}\label{1}
\Ord\inlineadj{F}{U}\Sup,
\end{equation}
has the remarkable feature that both functors involved are embeddings. This allows us to view $\Sup$ as a part of $\Ord$, but also $\Ord$ as a part of $\Sup$. The first viewpoint corresponds to the common understanding that a complete lattice is an ordered set in which all suprema exist and that a sup-morphism is an order-preserving function that preserves suprema. More technically: $\Sup$ is the category of Eilenberg-Moore algebras for the `free-cocompletion KZ-doctrine' on $\Ord$, which sends an ordered set to the set of its downclosed subsets ordered by inclusion. The second point of view is what the notion of totally algebraic complete lattice is all about. Recall that an element $a$ of a complete lattice $L$ is {\em totally compact} (a.k.a.\ {\em supercompact}) when, for any downclosed $A\subseteq L$, $a\leq\bigvee A$ implies $a\in A$; and a complete lattice $L$ is {\em totally algebraic} (a.k.a.\ {\em superalgebraic}) when each element is the supremum of totally compact ones [Gierz\etal, 1980]. It turns out that the replete image of the left adjoint in the above adjunction is precisely the subcategory of $\Sup$ of totally algebraic objects and left adjoint morphisms; thus $\Ord$ is described intrinsically (and up to equivalence) as a part of $\Sup$.
\par
We want to broaden the situation depicted above: instead of studying ordered sets, i.e.\ ordered sheaves on the two-element Boolean algebra $\2$, we want to consider ordered sheaves on any small quantaloid $\Q$. The latter form a 2-category $\Ord(\Q)$ that we define\footnote{See [Stubbe, 2005b] for a more ``elementary'' definition of ordered sheaves on a quantaloid.} as $\Cat\cc(\Q\si)$: its objects are the Cauchy-complete categories enriched in the split-idempotent completion of the quantaloid $\Q$, and its morphisms are the $\Q$-enriched functors. It was proved by I. Stubbe [2007b] that the category of ``internal sup-lattices and sup-morphisms'' in $\Ord(\Q)$ is (biequivalent to) $\Mod(\Q):=\QUANT(\Q\op,\Sup)$, the quantaloid of $\Q$-modules. That is to say, in analogy with the situation in \eqref{1} above, there is a biadjunction
\begin{equation}\label{2}
\Ord(\Q)\inlineadj{F}{U}\Mod(\Q)
\end{equation}
that splits the free-cocompletion KZ-doctrine on $\Ord(\Q)$ and such that moreover $\Mod(\Q)$ is (biequivalent to) the category of Eilenberg-Moore algebras for that doctrine. 
This describes $\Mod(\Q)$ as part of $\Ord(\Q)$, and for $\Q=\2$ we thus recover exactly half of the situation described in the first paragraph above. But what about the other half: Can we also intrinsically characterise $\Ord(\Q)$ as part of $\Mod(\Q)$? Can we give a module-theoretic condition on an object of $\Mod(\Q)$ that makes it equivalent to the free cocompletion of an object of $\Ord(\Q)$? And how do morphisms then relate?
\par
With our Definition \ref{20} and our Theorem \ref{19} we answer these questions in the affirmative: we prove that $\Map(\Mod\lpg(\Q))$, defined as the subcategory of $\Mod(\Q)$ of {\em locally principally generated $\Q$-modules and left adjoint $\Q$-module morphisms}, is precisely the replete image of the left biadjoint functor in \eqref{2} above:
$$\xymatrix@R=15mm@C=12mm{
\Ord(\Q)\ar[r]^F\ar@{=}[dr] & \Mod(\Q) \\
 & \Map(\Mod\lpg(\Q))\ar@{^{(}->}[u]}$$
\par
The technology that we use to solve this problem is $\Q$-enriched categorical algebra, as pioneered (in greater generality) by [B\'enabou, 1967; Walters, 1981; Street, 1983a] and more recently surveyed by [Stubbe, 2005a]. More particularly, in this paper we build further on results from [Stubbe, 2007a], which treats totally algebraic cocomplete $\Q$-categories, and [Stubbe, 2006], where an explicit comparison is given between cocomplete $\Q$-categories and $\Q$-modules.

\subsection*{Skew local homeomorphisms}

The notion of ordered sheaf on a quantaloid $\Q$ is so devised that, when taking $\Q$ to be the one-object suspension of a locale $X$ (i.e.\ $\Q$ has one object, the elements of $X$ are viewed as arrows of $\Q$, composition of which corresponds to finite meets in $X$, the identity arrow thus being the top element $\top$ of $X$), $\Ord(X)$ is equivalent to the category of ordered objects and order-preserving morphisms in the topos $\Sh(X)$ of sheaves on $X$ [Walters, 1981; Borceux and Cruciani, 1998; Stubbe, 2005b]: $\Ord(X)\simeq\Ord(\Sh(X))$. Our general results on $\Q$-modules from the first part of this paper surely specialise to the localic case: ordered sheaves on $X$ can be described equivalently as locally principally generated $X$-modules, order-preserving morphisms then correspond to left adjoint $X$-module morphisms.
\par
In analogy with ring theory it is very natural to regard a locale morphism $f\: Y\to X$ as a left $X$-module $(Y,\circ_f)$ with action $y\circ_f x:=y\wedge f^*(x)$ [Joyal and Tierney, 1984]. This construction extends to a (contravariant) embedding of the slice category $\Loc\slice{X}$ in $\Mod(X)$. Thus it is natural to try to characterise the subcategory of $\Loc\slice{X}$ which corresponds, under this embedding, to the locally principally generated $X$-modules and the left adjoint $X$-module morphisms between them, or in other words, to $\Ord(X)$. In Definitions \ref{S6} and \ref{S17} we introduce the locale theoretic notions of {\em skew open morphism} and {\em skew local homeomorphism}, and in Theorem \ref{S20} we then prove that $(\Loc\slice{X})\o\slh$, by definition the (non-full) subcategory of $\Loc\slice{X}$ of skew local homeomorphisms as objects and skew open morphisms between them, is the sought-after equivalent of $\Ord(X)$.
\par
A local homeomorphism is necessarily a skew local homeomorphism; and an open locale morphism is always skew open too. Thus the category $\LH\slice{X}$ of local homeomorphisms over $X$ is naturally a full subcategory of $(\Loc\slice{X})\o\slh$. This situation too can be stated in terms of $X$-modules: in Definition \ref{S22} we introduce {\em \'etale $X$-modules} as a particular kind of locally principally generated $X$-modules, such that in Theorem \ref{S23} we can prove that the full subcategory $\Map(\Mod\et(X))$ of $\Map(\Mod\lpg(X))$ defined by the \'etale $X$-modules is indeed equivalent to $\LH\slice{X}$.
\par
The category $\LH\slice{X}$ is a well-known equivalent of the topos $\Sh(X)$, see e.g.\ [Mac~Lane and Moerdijk, 1992, p.~524]; thus the inclusion of local homeomorphisms over $X$ into skew local homeomorphisms over $X$, or equivalently the inclusion of \'etale $X$-modules into locally principally generated $X$-modules, is precisely the same thing as the inclusion of sheaves on $X$ into ordered sheaves on $X$. Moreover, in the realm of enriched categorical structures it is a matter of fact that the ``$X$-sets'' of M. Fourman and D. Scott [1979] (see also [Borceux, 1994, Vol.~3; Borceux and Cruciani, 1998; Johnstone, 2002, p.~502--513]) are included in $X$-orders. All this establishes the following unifying diagram of equivalent embeddings of categories of symmetric (or discrete) objects into 2-categories of asymmetric (or ordered) objects:
$$\xymatrix@=12mm{
\Ord(\Sh(X))\ar@{=}[r] & \Ord(X)\ar@{=}[r] & (\Loc\slice{X})\o\slh\ar@{=}[r] & \Map(\Mod\lpg(X)) \\
\Sh(X)\ar@{=}[r]\ar@{^{(}->}[u] & \Set(X)\ar@{=}[r]\ar@{^{(}->}[u] & \LH\slice{X}\ar@{=}[r]\ar@{^{(}->}[u] & \Map(\Mod\et(X))\ar@{^{(}->}[u]}$$
This shows the relation between (ordered) sheaves on a locale $X$ as (i) functors on $X$ satisfying gluing conditions, (ii) $X$-enriched categorical structures, (iii) locale morphisms into $X$ and (iv) $X$-modules.

\subsection*{Overview of contents}

Sections \ref{B} through \ref{F} of this paper are concerned with the translation of the notion of {\em orderd sheaf on a small quantaloid $\Q$} from its original definition in terms of $\Q$-enriched categorical structures [Walters, 1981; Stubbe, 2005b] to the language of $\Q$-modules. To make this paper self-contained, we therefore start with an overview of the $\Q$-enriched categorical algebra that we need: in Section \ref{B} we recall the definition of $\Q$-categories, functors and distributors; we speak of weighted colimits in a $\Q$-category and of Cauchy-completeness of a $\Q$-category; and we end with the definition of $\Q$-order. We have tried to include the relevant ``historical'' references, but in practice we refer mostly to the more recent [Stubbe, 2005a, 2005b, 2006] whose notations we follow. In Section \ref{C} we recall some material on totally algebraic cocomplete $\Q$-categories from [Stubbe, 2007a]; in fact, we recast the definition of a totally compact object in a way that suits our needs further on. (In Section \ref{AA} we explain a biadjunction involving totally algebraic cocomplete $\Q$-categories: strictly speaking it is of no technical importance for the rest of this paper, but since it may be of independent interest we have chosen to add it as an Addendum.) Section \ref{D} then contains the crucial translation from $\Q$-enrichment to $\Q$-variation -- to borrow a term from [Betti {\it et al.}, 1983], later picked up by [Gordon and Power, 1997] and [Stubbe, 2006] -- where we introduce the notion of {\em principally generated $\Q$-module}. Our first main theorem is that Cauchy-complete $\Q$-categories and functors between them form a category which is equivalent to that of principally generated $\Q$-modules and left adjoint module morphisms. Finally, in Section \ref{E} we explain how a so-called {\em locally principally generated $\Q$-module} is the same thing as a principally generated module on the split-idempotent completion of $\Q$, thus paving the way for our second main result: $\Q$-orders (meaning Cauchy-complete categories enriched in the split-idempotent completion of $\Q$) and their morphisms are essentially the same thing as locally principally generated $\Q$-modules and left adjoint module morphisms. We discuss some examples in Section \ref{F}.
\par
The second part of this paper, contained in Section \ref{SX}, is devoted to the application of the above to the specific case where $\Q$ is the one-object suspension of a locale $X$ (viewed as monoid $(X,\wedge,\top)$). Locally principally generated $X$-modules are then equivalent to ordered sheaves on $X$, which this time can really be understood as ordered objects in the topos $\Sh(X)$ [Walters, 1981; Borceux and Cruciani, 1998; Stubbe, 2005b]. (If one takes for granted that locally principally generated $X$-modules are ordered sheaves on $X$, then one can start reading Section \ref{SX} right away; this second part of the paper is technically speaking rather independent from the first part.) But locally principally generated $X$-modules can also be expressed in terms of certain locale morphisms into $X$, and it is their study that we deal with here. Thus, we begin by briefly explaining, taking hints from [Joyal and Tierney, 1984], how any locale morphism into $X$ can be regarded as an $X$-module; we define {\em skew open morphisms} in the slice category $\Loc\slice{X}$ to correspond to left adjoint $X$-module morphisms. A detailed study of locally principally generated $X$-modules is carried out thereafter; we show in particular that any such $X$-module is necessarily induced by a locale morphism into $X$. This work being done, we come to our third main result of the paper: we define {\em skew local homeomorphisms} in terms of coverings by {\em skew open sections}, and prove that the subcategory of $\Loc\slice{X}$ with skew local homeomorphisms as objects and skew open morphisms between them, is equivalent to the category of locally principally generated $X$-modules and left adjoint module morphisms, {\it viz.}\ ordered sheaves on $X$. Remarking that local homeomorphisms are necessarily skew local homeomorphisms, we end with the identification of {\em \'etale $X$-modules} as those locally principally generated $X$-modules which correspond to local homeomorphisms, {\it viz.}\ sheaves on $X$.

\subsection*{Further work}

In this paper we do not speak of the ``internal logic'' of the category of ordered sheaves on a quantaloid -- it is still an open problem -- but we hope that our contribution here will be helpful to investigate this. In this respect it should be interesting to investigate links with examples of noncommutative logics developed from a more logical (rather than algebraic, i.e.\ sheaf theoretic) point of view: for example, R. Goldblatt's [2006] encoding of predicates in some non-commutative logic as quantale-valued functions on a set (which can be seen as elements of a principally generated module!); or [Baltag {\it et al.}, 2007], who use quantale-modules in their treatment of epistemic logic; or K. Rosenthal's [1994] model for the ``bang'' operator in linear logic via modules on a quantale; the construction in [Coniglio and Miraglia, 2001] of a logic from a very particular notion of sheaf on a restricted class of quantales; or the quantale based semantics for propositional normal modal logic in [Marcelino and Resende, 2008]; and many others.

In [Resende and Rodrigues, 2008] local homeomorphisms into $X$ are shown to correspond to Hilbert $X$-modules with a Hilbert basis (a special case of Hilbert $Q$-modules, for $Q$ an involutive quantale [Paseka, 1999]). Their results and our results in Section \ref{SX} (particularly Definition \ref{S22}) are similar in that we both provide a description of local homeomorphisms into a locale $X$ in terms of particular $X$-modules. In [Heymans and Stubbe, 2008] we explain, in the generality of modules on an involutive quantale $Q$, the precise relationship between Hilbert $Q$-modules admitting a Hilbert basis on the one hand, and locally principally generated $Q$-module satisfying a suitable {\em symmetry condition} on its locally principal elements on the other hand; we argue that the latter are precisely sets with a $Q$-valued {\em equality}, i.e.\ a ``$Q$-sets'' rather than a ``$Q$-orders''. Applied to a locale $X$ this gives (ordinary) sheaves on $X$; better still, applied to suitably constructed involutive quantales we can describe all Grothendieck toposes, in a manner closely related to [Walters, 1982].
 
Our current work continues along this line and focuses on $\Q$-valued equalities in the generality of an involutive quantaloid $\Q$, more specifically on the interplay between symmetric and non-symmetric $\Q$-categories and their Cauchy completions. This is directly related to [Walters, 1982; Betti and Walters, 1982; Freyd and Scedrov, 1990] but also has ties with [Gylys, 2001; Heymans, 2009].

\section{Preliminaries}\label{B}

\subsection*{Quantaloids}

Let $\Sup$ denote the category of complete lattices and maps that preserve arbitrary suprema ({\em suplattices} and {\em supmorphisms}): it is symmetric monoidal closed for the usual tensor product. A {\em quantaloid} is a $\Sup$-enriched category; a one-object quantaloid is most often thought of as a monoid in $\Sup$: it is a {\em quantale}. A $\Sup$-functor between quantaloids is sometimes called a {\em homomorphism}; $\QUANT$ denotes the (``illegitimate'') category of quantaloids and their homomorphisms. The standard reference on categories enriched in a symmetric monoidal category in general is [Kelly, 1982]; for quantales and quantaloids in particular there is [Rosenthal, 1990, 1996].
\par
Composition with a morphism $f\:X\to Y$ in a quantaloid $\Q$ gives rise to adjunctions, one for each $A\in\Q$,
$$\Q(A,X)\xymatrix@=15mm{\ar@{}[r]|{\perp}\ar@<1mm>@/^2mm/[r]^{f\circ-}&\ar@<1mm>@/^2mm/[l]^{[f,-]}}\Q(A,Y)\hspace{3ex}\mbox{ 
and }\hspace{3ex} 
\Q(Y,A)\xymatrix@=15mm{\ar@{}[r]|{\perp}\ar@<1mm>@/^2mm/[r]^{-\circ 
f}&\ar@<1mm>@/^2mm/[l]^{\{f,-\}}}\Q(X,A).$$
These right adjoints are respectively called {\em lifting} and {\em 
extension} (through $f$). We shall keep the notations ``$[-,-]$'' and ``$\{-,-\}$'' for liftings and extensions in any quantaloid that follows; no confusion shall arise\footnote{These right adjoints also go by the name of {\em residuations} when $\Q$ is a quantale, i.e.\ a one-object quantaloid. Whereas our notations are the usual ones in category theory (for closed monoidal categories or bicategories), other notations instead of $[f,g]$ and $\{f,g\}$ that can be found in the literature include $f\rightarrow g$ and $g\leftarrow f$, or $f\rightarrow_r g$ and $f\rightarrow_l g$, or $f/g$ and $f\backslash g$, or $g\raiseletter{f}$ and $\raiseletter{f}g$.}.
\par
Given a quantaloid $\Q$ we write $\Q\si$ (``{\sf si}'' stands for ``split the idempotents'') for the new quantaloid whose objects are the idempotent arrows in $\Q$, and in which an arrow from an idempotent $e\:A\to A$ to an idempotent $f\:B\to B$ is a $\Q$-arrow $g\:A\to B$ satisfying $g\circ e=g=f\circ g$. Composition in $\Q\si$ is done as in $\Q$, the identity in $\Q\si$ on some idempotent $e\:A\to A$ is $e$ itself, and the local order in $\Q\si$ is that of $\Q$. There is an obvious inclusion $j\:\Q\to\Q\si$, mapping $f\:A\to B$ to $f\:1_A\to 1_B$, which expresses $\Q\si$ as the universal split-idempotent completion of $\Q$ in $\QUANT$.
\par
When $\Q$ is a small quantaloid, $\Mod(\Q)$ shall be shorthand for $\QUANT(\Q\op,\Sup)$: the objects of this (large) quantaloid are called the {\em modules on $\Q$}, or briefly $\Q$-modules. Since idempotents split in $\Sup$, and noting that $\Q\si$ is small whenever $\Q$ is, it follows that $\Mod(\Q)$ is equivalent to $\Mod(\Q\si)$. We shall come back to modules on $\Q$ and on $\Q\si$ in Sections \ref{D} and \ref{E}.
\par
A quantaloid is in particular a locally ordered category, and therefore we can straightforwardly define {\em adjoint pairs} in any quantaloid: $f\:X\to Y$ is left adjoint to $g\:Y\to X$ (and $g$ is right adjoint to $f$, written $f\dashv g$) when $1_X\leq g\circ f$ and $f\circ g\leq 1_X$. If a morphism $f\:X\to Y$ has a right adjoint, then the latter is unique, and we shall often use $f^*$ as its notation. Because left adjoints are sometimes called ``maps'', $\Map(\Q)$ is our notation for the (locally ordered) category of left adjoints in a quantaloid $\Q$.

\subsection*{Quantaloid-enriched categories}

A quantaloid is a bicategory and therefore it may serve itself as base for enrichment [B\'enabou 1967; Walters, 1981; Street, 1983a]. The theory of quantaloid-enriched categories, functors and distributors is surveyed in [Stubbe, 2005a] where many more appropriate references are given; here we can only provide a brief summary, but we follow the notations of {\it op.\ cit.}\ for easy cross reference. To avoid size issues we work from now on {\em with a small quantaloid $\Q$}. 
\par
A {\em $\Q$-category} $\bbA$ consists of a set $\bbA_0$ of `objects', a `type' function $t\:\bbA_0\to\Q_0$, and for any $a,a'\in\bbA_0$ a `hom-arrow' $\bbA(a',a)\:ta\to ta'$ in $\Q$; these data are required to satisfy 
$$\bbA(a'',a')\circ\bbA(a',a)\leq\bbA(a'',a)\mbox{\ \ \ and\ \ \ }1_{ta}\leq\bbA(a,a)$$ 
for all $a,a',a''\in\bbA_0$. A {\em functor} $F\:\bbA\to \bbB$ is a map $\bbA_0\to\bbB_0\:a\mapsto Fa$ that satisfies 
$$ta=t(Fa)\mbox{\ \ \ and\ \ \ }\bbA(a',a)\leq\bbB(Fa',Fa)$$ 
for all $a,a'\in\bbA_0$. 
\par
A $\Q$-category $\bbA$ has an {\em underlying order}\footnote{By an {\em order} we mean a reflexive and transitive relation, i.e.\ a (small) category with at most one arrow between any two objects; some call this a preorder. We shall speak of a {\em partial order} or an {\em antisymmetric order} if we require moreover antisymmetry.} $(\bbA_0,\leq)$: for $a,a'\in\bbA_0$ define $a\leq a'$ to mean $ta=ta'=:A$ and $1_A\leq\bbA(a,a')$. If $a\leq a'$ and $a'\leq a$ we write $a\cong a'$ and say that these are {\em isomorphic objects} in $\bbA$. For parallel functors $F,G\:\bbA\biar\bbB$ we now put $F\leq G$ when $Fa\leq Ga$ for every $a\in\bbA_0$. With the obvious composition and identities we thus obtain a locally ordered category $\Cat(\Q)$ of $\Q$-categories and functors. Precisely because $\Cat(\Q)$ is a 2-category, we can from now on unambiguously use 2-categorical notions such as adjoint functors, Kan extensions, and so on.
\par
To give a {\em distributor} (or {\em module} or {\em profunctor}) $\Phi\:\bbA\dist\bbB$ between $\Q$-categories is to specify, for all $a\in\bbA_0$ and $b\in\bbB_0$, arrows $\Phi(b,a)\:ta\to tb$ in $\Q$ such that 
$$\bbB(b,b')\circ\Phi(b',a)\leq\Phi(b,a)\mbox{\ \ \ and\ \ \ }\Phi(b,a')\circ\bbA(a',a)\leq\Phi(b,a)$$
for every $a,a'\in\bbA_0$, $b,b'\in\bbB_0$. Two distributors $\Phi\:\bbA\dist\bbB$, $\Psi\:\bbB\dist\bbC$ 
compose: we write $\Psi\tensor\Phi\:\bbA\dist\bbC$ for the 
distributor with elements
$$\Big(\Psi\tensor\Phi\Big)(c,a)=\bigvee_{b\in\bbB_0}\Psi(c,b)\circ\Phi(b,a).$$
The identity distributor on a $\Q$-category $\bbA$ is 
$\bbA\:\bbA\dist\bbA$ itself, i.e.\ the distributor with elements 
$\bbA(a',a)\:ta\to ta'$. For parallel distributors $\Phi,\Phi'\:\bbA\bidist\bbB$ we define $\Phi\leq\Phi'$ to mean that $\Phi(b,a)\leq\Phi'(b,a)$ for every $a\in\bbA_0$, $b\in\bbB_0$. It is easily seen that $\Q$-categories and distributors form a quantaloid $\Dist(\Q)$. 
\par
Every functor $F\:\bbA\to\bbB$ between $\Q$-categories represents an adjoint pair of distributors: 
\begin{itemize}
\itemsep=-4pt
\item the left adjoint $\bbB(-,F-)\:\bbA\dist\bbB$ 
has elements 
$\bbB(b,Fa)\:ta\to tb$,
\item the right adjoint $\bbB(F-,-)\:\bbB\dist\bbA$ 
has elements 
$\bbB(Fa,b)\:tb\to ta$.
\end{itemize}
The assignment $F\mapsto\bbB(-,F-)$ is a (bijective and) faithful 2-functor from 
$\Cat(\Q)$ to $\Dist(\Q)$; thus, whenever a distributor 
$\Phi\:\bbA\dist\bbB$ is represented by a functor $F\:\bbA\to\bbB$, 
this $F$ is essentially unique.

\subsection*{Weighted colimits}

In a $\Q$-enriched category $\bbC$ we can speak of weighted limits and colimits, as introduced by R. Street [1983a] for general bicategory-enrichment. For our short account here we use [Stubbe, 2005a, 2006] as references, but we should certainly also mention the work of [Gordon and Power, 1997] on conical (co)limits and (co)tensors (to be explained below).
\par
Given a distributor $\Phi\:\bbA\dist\bbB$ and a functor $F\:\bbB\to\bbC$, a functor $K\:\bbA\to\bbC$ is the {\em $\Phi$-weighted colimit of $F$} when it satisfies
$$\bbC(K-,-)=[\Phi,\bbC(F-,-)]$$
in $\Dist(\Q)$ (and in that case it is essentially unique, so we write it as $\colim(\Phi,F)$). The $\Q$-category $\bbC$ is {\em cocomplete} when it admits all such weighted colimits. A functor $H\:\bbC\to\bbC'$ is {\em cocontinuous} when it preserves all colimits that happen to exist in $\bbC$: 
$H\circ\colim(\Phi,F)\cong\colim(\Phi,H\circ F)$. A left adjoint 
functor is always cocontinuous; conversely, if the domain of a 
cocontinuous functor is cocomplete, then that functor is left 
adjoint. Cocomplete $\Q$-categories and cocontinuous functors form a 
sub-2-category $\Cocont(\Q)$ of $\Cat(\Q)$. (Dually one can speak of {\em weighted limits}, {\em complete $\Q$-categories} and {\em continuous functors}; we shall not explicitly need these further on, however, it is a matter of fact that a $\Q$-category is complete if and 
only if it is cocomplete [Stubbe, 2005a, 5.10].)
\par
Every object $X$ of a quantaloid $\Q$ determines a one-object 
$\Q$-category $*_X$ whose single hom-arrow is $1_X$. A {\em 
contravariant presheaf of type $X$} on a $\Q$-category $\bbA$ is a 
distributor $\phi\:*_X\dist\bbA$; these are the objects of a 
cocomplete $\Q$-category $\P\bbA$ whose hom-arrows are given by 
lifting in $\Dist(\Q)$. Every object $a\in\bbA_0$ determines, and is 
determined by, a functor $\Delta a\:*_{ta}\to\bbA$; thus $a\in\bbA_0$ also 
represents a (left adjoint) presheaf $\bbA(-,a)\:*_{ta}\dist\bbA$. 
The {\em Yoneda embedding} 
$Y_{\bbA}\:\bbA\to\P\bbA\:a\mapsto\bbA(-,a)$ is a fully 
faithful\footnote{A functor $F\:\bbA\to\bbB$ is {\em fully faithful} 
when $\bbA(a',a)=\bbB(Fa',Fa)$ for every $a,a'\in\bbA_0$.} 
continuous functor. The presheaf construction $\bbA\mapsto\P\bbA$ 
extends to a 2-functor $\P\:\Cat(\Q)\to\Cocont(\Q)$ which is left 
biadjoint to the inclusion 2-functor, with the Yoneda embeddings as 
unit; thus presheaf categories are the freely cocomplete ones. In fact, a $\Q$-category $\bbC$ is cocomplete if and only if the Yoneda embedding admits a left adjoint in $\Cat(\Q)$; if this is the case we write $\sup_{\bbC}\:\P\bbC\to\bbC$ for that left adjoint: it maps a presheaf $\phi$ on $\bbC$ to the weighted colimit $\sup_{\bbC}(\phi):=\colim(\phi,1_{\bbC})$. Note by the way that $\sup_{\bbC}\circ Y_{\bbC}\cong 1_{\bbC}$; actually, $Y_{\bbC}$ admits a left adjoint if and only if it admits a left inverse.
\par
Given any $X\in\Q$ we shall write $(\bbC_X,\leq)$ for the (possibly empty) sub-order of $(\bbC_0,\leq)$ containing all $c\in\bbC_0$ for which $tc=X$. If for each $X\in\Q$ the order $(\bbC_X,\leq)$ is a complete lattice, we say that $\bbC$ is {\em order-cocomplete}. On the other hand, for each morphism $f\:A\to B$ in $\Q$ we can consider the one-element distributor $(f)\:*_A\dist*_B$. Suppose that $c\in\bbC$ is of type $tc=B$, and that $\colim((f),\Delta c)$ exists: it is itself a functor from $*_A$ to $\bbC$, and can thus be identified with an element of $\bbC$ of type $A$. We write that element as $c\tensor f$ and call it the {\em tensor} of $c$ with $f$; $\bbC$ {\em has all tensors} when all such colimits with one-element weights exist. The dual notion is {\em cotensor}. It has been proved in [Stubbe, 2006, 2.13] that a $\Q$-category $\bbC$ is cocomplete if and only if it is order-cocomplete and has all tensors and cotensors; moreover, for any $\Phi\:\bbA\dist\bbB$ and $F\:\bbB\to\bbC$, the weighted colimit $\colim(\Phi,F)\:\bbA\to\bbC$ is the functor defined by $a\mapsto\bigvee_{b\in\bbB_0}Fb\tensor\Phi(b,a)$.

\subsection*{Ordered sheaves on a quantaloid}

The importance of the following notion has first been recognised by B. Lawvere [1973] in the context of categories enriched in a monoidal category.
\par
A $\Q$-category $\bbC$ is {\em Cauchy complete} if for any other $\Q$-category $\bbA$ the map
$$\Cat(\Q)(\bbA,\bbC)\to\Map(\Dist(\Q))(\bbA,\bbC)\:F\mapsto \bbC(-,F-)$$
is surjective, i.e.\ when any left adjoint distributor (also called {\em Cauchy distributor}) into $\bbC$ is represented by a functor. This is equivalent to the requirement that $\bbC$ admits any colimit weighted by a Cauchy distributor; and moreover such weighted colimits are {\em absolute} in the sense that they are preserved by any functor [Street, 1983b]. We write $\Cat\cc(\Q)$ for the full subcategory of $\Cat(\Q)$ whose objects are the Cauchy complete $\Q$-categories.
\par
Now we have everything ready to define the central notion of this paper [Stubbe, 2005b].
\begin{definition}
For a small quantaloid $\Q$, we write $\Ord(\Q)$ for the locally ordered category $\Cat\cc(\Q\si)$, and call its objects {\em ordered sheaves on $\Q$}, or simply {\em $\Q$-orders}.
\end{definition}
When taking $\Q$ to be the (one-object suspension of) a locale $X$, $\Ord(X)$ is the category of ordered objects and order-preserving morphisms in the topos $\Sh(X)$, as first (implicitly) observed by B. Walters [1981] (but see also [Borceux and Cruciani, 1998] for the locale-specific notion, and [Stubbe, 2005b] for the generalisation to quantaloids and the comparison between [Walters, 1981] and [Borceux and Cruciani, 1998]). Obviously, this example inspired our terminology.

\section{Total algebraicity revisited}\label{C}

We shall review and expand the material that we need from [Stubbe, 2007a].
\begin{definition}[Stubbe, 2007a]\label{3} Let $\bbA$ be a cocomplete $\Q$-category. The {\em totally below distributor} on $\bbA$ is the right extension of $\bbA(-,\sup_{\bbA}-)$ through $\P\bbA(Y_{\bbA}-,-)$ in $\Dist(\Q)$:
$$\xy\xymatrix@=20mm{
\P\bbA\ar[r]|{\distsign}^{\P\bbA(Y_{\bbA}-,-)}\ar[d]|{\distsign}_{\bbA(-,\sup_{\bbA}-)} & \bbA \\
\bbA\ar@{.>}[ur]|{\distsign}_{\Theta_{\bbA}:=\Big\{\bbA(-,\sup_{\bbA}-),\P\bbA(Y_{\bbA}-,-)\Big\}.}}\POS(8,-10)\drop{\leq}\endxy$$
An object $a\in\bbA$ is {\em totally compact} when $1_{ta}\leq\Theta_{\bbA}(a,a)$. Writing $i_{\bbA}\:\bbA\c\to\bbA$ for the full embedding of the totally compacts, $\bbA$ is {\em totally algebraic} when the left Kan extension of $i_{\bbA}$ along itself is isomorphic to $1_{\bbA}$.
\end{definition}
In the simplest possible case, when $\Q$ is the (one-object suspension of) the two-element Boolean algebra $(\2,\wedge,\top)$, a $\Q$-category $\bbA$ is an ordered set $(A,\leq)$, it is cocomplete precisely when $(A,\leq)$ is a sup-lattice, and the distributor $\Theta_{\bbA}$ is the following ``totally below'' relation: $a'\lll a$ when, for every down-closed subset $D\subseteq A$, $a\leq\bigvee D$ implies $a'\in D$. A totally compact element is one which is totally below itself, and $(A,\leq)$ is totally algebraic if and only if every element is the supremum of the totally compacts below it. These notions are related to, but stronger than, the ``way below'' relation and the ``algebraic'' sup-lattices [Gierz {\it et al.}, 1980].
\begin{theorem}[Stubbe, 2007a]\label{3.1}
The 2-functor
$$\P\:\Dist(\Q)\to\Cocont(\Q)\:\Big(\Phi\:\bbA\dist\bbB\Big)\mapsto
\Big(\Phi\tensor-\:\P(\bbA)\to\P(\bbB)\Big)$$
is locally an equivalence, and its corestriction to the full sub-2-category of totally algebraic cocomplete $\Q$-categories is a biequivalence: $\Dist(\Q)\simeq\Cocont\ta(\Q)$.
\end{theorem}
We may restrict the biequivalence of which this theorem speaks, to left adjoints: we then obtain the biequivalence (which we write with the same letter)
$$\P\:\Map(\Dist(\Q))\to\Map(\Cocont\ta(\Q)).$$
But the definition of Cauchy completeness for $\Q$-categories implies that also
$$\Cat\cc(\Q)\to\Map(\Dist(\Q))\:\Big(F\:\bbA\to\bbB\Big)\mapsto
\Big(\bbB(-,F-)\:\bbA\dist\bbB\Big)$$
is a biequivalence, hence composing these two we get a third biequivalence (which we still write with the same letter):
\begin{corollary}\label{3.2}
The 2-functor
$$\P\:\Cat\cc(\Q)\to\Map(\Cocont\ta(\Q))\:\Big(F\:\bbA\to\bbB\Big)\mapsto
\Big(\bbB(-,F-)\tensor-\:\bbA\dist\bbB\Big)$$
is a biequivalence.
\end{corollary}
The inverse biequivalence is given by ``taking totally compact objects''. More precisely, if $F\:\bbA\to\bbB$ is in $\Map(\Cocont(\Q))$ then it maps totally compact objects of $\bbA$ to totally compact objects of $\bbB$, hence we get a functor $F\c\:\bbA\c\to\bbB\c$ out of it. If $\bbA$ and $\bbB$ are moreover totally algebraic, then $\bbA\c$ and $\bbB\c$ are Cauchy complete. This describes a 2-functor $(-)\c\:\Map(\Cocont\ta(\Q))\to\Cat\cc(\Q)$, which turns out to be the sought-after inverse.
\par
In fact, the biequivalence in Corollary \ref{3.2} can also be seen as resulting from (co)restricting the following biadjunction to the objects for which the (co)unit is an equivalence: 
$$\Cat(\Q)\inlineadj{\P}{\ (-)\c}\Map(\Cocont(\Q)).$$
This is without importance for what follows, so we shall not include the details here; but since this may be of independent interest, we have written the details in a technical Addendum at the end of this paper.
\par
Several equivalent expressions for the definition of totally compact object are given in [Stubbe, 2007a] but for the purposes of this paper the following are particularly useful\footnote{Especially condition \eqref{4.3} in Proposition \ref{4} is reminiscent of the notion of {\em atom} defined by M. Bunge [1969] and that of {\em small-projective object} defined by M. Kelly [1982].}:
\begin{proposition}\label{4}
Let $a\in\bbA$ be an object, of type $A\in\Q$ say, of a cocomplete $\Q$-category. The following conditions are equivalent:
\begin{enumerate}
\item\label{4.1} $a$ is totally compact,
\item\label{4.2} for all $\phi\in\P\bbA$: $\phi(a)=\bbA(a,\sup_{\bbA}(\phi))$,
\item\label{4.3} the functor $H_a\:\bbA\to\P(*_A)\:x\mapsto\bbA(a,x)$ (``homming with $a$'') is cocontinuous,
\item\label{4.4} the functor $T_a\:\P(*_A)\to\bbA\:f\mapsto a\tensor f$ (``tensoring with $a$'') is cocontinuous and admits a cocontinuous right adjoint.
\end{enumerate}
\end{proposition}
\proof
(\ref{4.1}$\iff$\ref{4.2}) It is easily seen (and spelled out in [Stubbe, 2007a, {\bf 5.2}]) that, for any $x,y\in\bbA$,
$$\Theta_{\bbA}(x,y)=\bigwedge_{\phi\in\P\bbA}\left\{\bbA(y,\sup_{\bbA}(\phi)),\phi(x)\right\}$$
and hence, straightforwardly,
\begin{eqnarray*}
1_{A}\leq\Theta_{\bbA}(a,a)
 & \iff & \forall\phi\in\P\bbA: 1_{A}\leq\left\{\bbA(a,\sup_{\bbA}(\phi)),\phi(a)\right\} \\
 & \iff & \forall\phi\in\P\bbA: \bbA(a,\sup_{\bbA}(\phi))\leq\phi(a).
\end{eqnarray*}
But because $\phi\leq Y_{\bbA}(\sup_{\bbA}(\phi))$ is automatic, and thus $\phi(a)\leq\bbA(a,\sup_{\bbA}(\phi))$ always holds, $a$ being totally compact is indeed equivalent to the clause in statement (\ref{4.2}).
\par
(\ref{4.2}$\iff$\ref{4.3}) By a straightforward calculation (e.g. using [Stubbe, 2006, Corollary 2.15]) it is easily seen that, for any $\phi\in\P(\bbA)$, $\phi(a)$ is the $\phi$-weighted colimit of $H_a$: 
$$\colim(\phi,H_a)=\bigvee_{x\in\bbA}H_a(x)\tensor\phi(x)=\bbA(a,-)\tensor\phi=\phi(a).$$
On the other hand it is clear by definition of $H_a$ that $\bbA(a,\sup_{\bbA}(\phi))=H_a(\sup_{\bbA}(\phi))$. Thus the formula in the second statement of the lemma can be rewritten as
$$\mbox{for all }\phi\in\P\bbA\mbox{: }\colim(\phi,H_a)=H_a(\sup_{\bbA}(\phi)).$$
As follows straightforwardly from [Stubbe, 2005a, 5.4], this in turn is equivalent to $H_a$ preserving {\em all} weighted colimits.
\par
(\ref{4.3}$\iff$\ref{4.4}) Due to $\bbA$'s cocompleteness, we surely have that all tensors with the object $a$ exist in $\bbA$; thus (the dual of) Proposition 3.2 in [Stubbe, 2006] says that $H_a$ is the right adjoint to $T_a$ in $\Cat(\Q)$. (A direct verification is very easy too.) Because the cocontinuous functors between cocomplete $\Q$-categories are precisely the left adjoint ones in $\Cat(\Q)$ [Stubbe, 2005a, Proposition 6.8], the result follows directly.
\endofproof

\section{Principally generated modules}\label{D}

A $\Q$-module $\F\:\Q\op\to\Sup$ determines a cocomplete $\Q$-category $\bbA_\F$ whose set of objects is $\coprod_X\F X$, with types given by $tx=X$ if and only if $x\in\F X$, and hom-arrows 
\begin{equation}\label{6}
\bbA_\F(y,x)=\bigvee\{f\:tx\to ty\mid \F(f)(y)\leq x\}.
\end{equation}
Similarly, a module morphism $\alpha\:\F\tto\G$ determines a cocontinuous functor $$F_{\alpha}\:\bbA_{\F}\to\bbA_{\G}\:x\mapsto\alpha_{tx}(x).$$
This sets up the biequivalence of $\Mod(\Q)$ with $\Cocont(\Q)$, as studied in great detail in [Stubbe, 2006, Section 4] based on work by R. Gordon and J. Power [1997]. We wish to characterise, purely in terms of $\Q$-modules, those $\F\in\Mod(\Q)$ for which the corresponding $\bbA_{\F}\in\Cocont(\Q)$ is totally algebraic. In order to do so, we must introduce some notations.
\par
Let $\F\:\Q\op\to\Sup$ be a $\Q$-module, and suppose that $a\in\F(A)$ for some $A\in\Q$. We shall write 
$$\Q(-,A)\twocell[9mm]{\tau_a}\F$$
for the $\Sup$-natural transformation that $a\in\F(A)$ corresponds with by the $\Sup$-enriched Yoneda Lemma. Often we shall loosely speak of ``an element $a$ of $\F$'', and even write ``$a\in\F$'', where actually we should be more precise and stipulate that $a\in\F(A)$ for some $A\in\Q$.
\par
Because $\Mod(\Q)$ is a (large) quantaloid, we can compute extensions and liftings. It is straightforward to verify, with the aid of the $\Sup$-enriched Yoneda Lemma, that for a $\Q$-module $\F\:\Q\op\to\Sup$ and $x\in\F(X)$ and $y\in\F(Y)$, the right lifting of $\tau_x$ through $\tau_y$,
\begin{equation}\label{11.0}
\begin{array}{c}
\xy\xymatrix@C=15mm@R=5mm{
\Q(-,X)\ar@{=>}[rd]^{\tau_x}\ar@{:>}[dd]_{[\tau_y,\tau_x]} \\
 & \F \\
\Q(-,Y)\ar@{=>}[ru]_{\tau_y}}\POS(8,-10)\drop{\leq}\endxy
\end{array}
\end{equation}
is precisely represented by the $\Q$-morphism $\bbA_{\F}(y,x)$ in \eqref{6}.
\par
Of course it makes sense to speak of adjoints in $\Mod(\Q)$; if a $\Q$-module morphism $\alpha\:\F\tto\G$ has a right adjoint\footnote{Whether $\alpha$ has a right adjoint in $\Mod(\Q)$ or not, each of its components certainly has a right adjoint in $\Ord$, say $\alpha_X'\:\G(X)\to\F(X)$, and these always form a {\em lax} natural $\Ord$-transformation $\alpha'\:\G\tto\F$. If $\alpha$ has a right adjoint in $\Mod(\Q)$, then -- for reasons of unicity of adjoints in $\Ord$ -- it must be $\alpha'$. In other words, $\alpha$ has a right adjoint in $\Mod(\Q)$ if and only if the lax natural $\Ord$-transformation $\alpha'$ is strictly natural and its components preserve suprema.}, we shall usually denote it as $\alpha^*\:\G\tto\F$.
\begin{proposition}\label{11}
Let $\F\:\Q\op\to\Sup$ be a $\Q$-module, and $\bbA_{\F}$ the associated cocomplete $\Q$-cat\-e\-go\-ry. We have the following:
\begin{enumerate}
\item\label{11.a} $a\in\bbA_{\F}$ is totally compact if and only if $a\in\F$ and $\tau_a$ is a left adjoint.
\item\label{11.b} $\bbA_{\F}$ is totally algebraic if and only if 
\begin{equation}\label{11.1}
\mbox{for each $x\in\F$, }\tau_x=\bigvee\Big\{\tau_a\circ[\tau_a,\tau_x]\bigmid\mbox{$a\in\F$ and $\tau_a$ is a left adjoint}\Big\}.
\end{equation}
\end{enumerate}
\end{proposition}
\proof
(\ref{11.a}) A representable module $\Q(-,A)\:\Q\op\to\Sup$ corresponds under the biequivalence $\Mod(\Q)\simeq\Cocont(\Q)$ with $\P(*_A)$, and a $\Sup$-natural transformation $\tau_a\:\Q(-,A)\tto\F$ corresponds to the ``tensoring with $a$'' cocontinuous functor $T_a\:\P(*_A)\to\bbA_{\F}\:f\mapsto a\tensor f=F(f)(a)$. Thus the fourth statement in Proposition \ref{4} proves the claim made here. 
\par
(\ref{11.b}) With suitable application of the $\Sup$-enriched Yoneda Lemma, it is easily deduced from \eqref{11.a} that the formula in \eqref{11.1}, which is stated in terms of the $\Q$-module $\F$, says precisely the same thing as
$$\mbox{for all $x\in\bbA_{\F}$, }x=\bigvee\{a\tensor\bbA_{\F}(a,x)\mid a\in(\bbA_{\F})\c\},$$
which is stated in terms of the cocomplete $\Q$-category $\bbA_{\F}$. But the right hand side in this latter formula is the explicit way of writing the value in $x$ of the (pointwise) left Kan extension of $i_{\bbA_{\F}}\:(\bbA_{\F})\c\to\bbA_{\F}$ along itself (see e.g.\ [Stubbe, 2005a, p.\ 26] combined with [Stubbe, 2006, Corollary 2.15]): hence it says that this left Kan extension is the identity functor on $\bbA_{\F}$. 
\endofproof
The preceding result promps the following definition.
\begin{definition}\label{12}
Let $\F\:\Q\op\to\Sup$ be a $\Q$-module. An $a\in\F$ is a {\em principal element} if $\tau_a$ is a left adjoint in $\Mod(\Q)$. Writing the set of principal elements of $\F$ as $\F\pr$, $\F$ is {\em principally generated} if
\begin{equation}\label{12.1}
\mbox{for each $x\in\F$, }\tau_x=\bigvee\Big\{\tau_a\circ[\tau_a,\tau_x]\bigmid a\in\F\pr\Big\}.
\end{equation}
\end{definition}
We can add the following succinct characterisation.
\begin{proposition}
Let $\F\:\Q\op\to\Sup$ be a $\Q$-module. The set of principal elements of $F$ is
$$\F\pr=\Big\{\tau(1_X)\mid \mbox{$X\in\Q$ and $\tau\:\Q(-,X)\tto\F$ is a left adjoint}\Big\},$$
and $\F$ is principally generated if and only if
\begin{equation}\label{11.2}
\id_{\F}=\bigvee\Big\{\tau\circ\tau^*\bigmid\mbox{$X\in\Q$ and $\tau\:\Q(-,X)\tto\F$ is a left adjoint}\Big\}.
\end{equation}
\end{proposition}
\proof
The first part of the Proposition is an application of the $\Sup$-enriched Yoneda Lemma.
\par
As for the second part, for any $x\in\F$ we can compute, by general calculation rules for liftings in a quantaloid\footnote{In any quantaloid, if $f\:A\to B$ is left adjoint to $f^*\:B\to A$, then for any $g\:C\to B$ we have $[f,g]=f^*\circ g$.}, that 
$$\left(\bigvee\Big\{\tau_a\circ\tau_a^*\bigmid a\in\F\pr\Big\}\right)\circ\tau_x 
=\bigvee\Big\{\tau_a\circ\tau_a^*\circ\tau_x\bigmid a\in\F\pr\Big\} \\
=\bigvee\Big\{\tau_a\circ[\tau_a,\tau_x]\bigmid a\in\F\pr\Big\}.$$
If we assume \eqref{11.2}, i.e.\ if the bracketed expression on the far left equals $\id_{\F}$, then clearly \eqref{12.1} follows. Conversely, assuming \eqref{12.1}, i.e.\ assuming that the far right expression equals $\tau_x=\id_{\F}\circ\tau_x$, and this for every $x\in\F$, implies -- because the representables are generators in $\Mod(\Q)$ -- that the bracketed expression on the far left must be equal to $\id_{\F}$.
\endofproof
From this we can now deduce an elegant characterisation of principally generated $\Q$-modules, entirely in terms of universal constructions\footnote{In any quantaloid, products coincide with coproducts, and are often called {\em direct sums}. An object $Y$ is an {\em adjoint retract} of an object $X$ when there exists a left adjoint $p\:X\to Y$ whose right adjoint $s\:Y\to X$ is its splitting ($p\circ s=1_Y$).} in the cocomplete quantaloid $\Mod(\Q)$ (we thank S. Lack for a stimulating discussion on this topic).
\begin{corollary}\label{11.3}
An $\F\in\Mod(\Q)$ is principally generated if and only if it is the adjoint retract of a direct sum of representable $\Q$-modules.
\end{corollary}
\proof
This holds by application of Lemma \ref{11.4} below to the family of all left adjoint $\Q$-module morphisms from representable $\Q$-modules to $\F$; compare with \eqref{12.1}.
\endofproof
\begin{lemma}\label{11.4}
Let $F$ be an object in a quantaloid with direct sums and consider a family of left adjoints with codomain $F$:
$$\Big\{F\inlinerevadj{f_i^*}{f_i}R_i\Big\}_{i\in I}.$$
Writing $f\:F\to\oplus_i R_i$ and $f^*\:\oplus_i R_i\to F$ for the unique factorisations of the $f_i$'s and the $f_i^*$'s we have that $f\circ f^*=\bigvee_i(f_i\circ f_i^*)$ and $f^*\circ f\geq 1_{\oplus_i R_i}$. Hence $1_F=\bigvee_i(f_i\circ f_i^*)$ if and only if $f$ and $f^*$ express $F$ as an adjoint retract of $\oplus_i R_i$.
\end{lemma}
\proof
Recall that, in any quantaloid, an object $\oplus_iR_i$ is the direct sum of a family of objects $(R_i)_{i\in I}$ if and only if there are morphisms $s_j\:R_j\to\oplus_iR_i$ and $p_j\:\oplus_iR_i\to R_j$, for all $j\in I$, such that $\bigvee_i(s_i\circ p_i)=1_{\oplus_iR_i}$ and $p_j\circ s_i=\delta_{i,j}$ (where $\delta_{i,j}\:X_i\to X_j$ is zero when $i\neq j$ and the identity otherwise). It follows that also $p_j\circ f^*=f_j^*$ and $f\circ s_j=f_j$ hold. To prove the equality, we compute that
$$f\circ f^*=f\circ(\bigvee_i(s_i\circ p_i))\circ f^*=\bigvee_i((f\circ s_i)\circ(p_i\circ f^*))=\bigvee_i(f_i\circ f_i^*).$$
The inequality is (because we are dealing with a 2-categorical (co)product) equivalent to requiring $p_j\circ 1_{\oplus_i R_i}\circ s_k\leq p_j\circ f^*\circ f\circ s_k$ for all $j,k\in I$. But the left hand side equals $\delta_{j,k}$ and the right hand side equals $f^*_j\circ f_k$, and since it is true that $1_{R_j}\leq f^*_j\circ f_j$, we are done.
\endofproof
\par
We shall write $\Mod\pg(\Q)$ for the full subquantaloid of $\Mod(\Q)$ determined by the principally generated $\Q$-modules, and thus obtain:
\begin{theorem}\label{13}
For any small quantaloid $\Q$, the biequivalence $\Mod(\Q)\simeq\Cocont(\Q)$ restricts to a biequivalence $\Mod\pg(\Q)\simeq\Cocont\ta(\Q)$.
\end{theorem}
Combined with the earlier observation in Corollary \ref{3.2} we get:
\begin{corollary}\label{13.1}
For any small quantaloid $\Q$, the locally ordered categories $\Cat\cc(\Q)$ and $\Map(\Mod\pg(\Q))$ are biequivalent.
\end{corollary}
\par
Although slightly off-topic, we find it important to remark that Corollary \ref{11.3} implies that the principally generated modules on $\Q$ form a {\em closed class of colimit weights} in the sense of [Albert and Kelly, 1988; Kelly and Schmitt, 2005]; in fact, this class is the {\em closure} of the class of (weights for) direct sums and adjoint retracts. The general theory explained for $\cal V$-enriched categories in the cited references implies that, for any small quantaloid $\Q$, $\Mod\pg(\Q)$ is precisely the free cocompletion of $\Q$ for direct sums and adjoint retracts, or equivalently, the free cocompletion of $\Q$ for all colimits weighed by a principally generated module. Since we know that $\Dist(\Q)\simeq\Cocont\ta(\Q)\simeq\Mod\pg(\Q)$, this at once describes the universal property of the distributor quantaloid too. In [Stubbe, 2005a] it is shown that $\Dist(\Q)$ is the universal ``direct sum and split monad'' completion of $\Q$; but it is trivial that, in a quantaloid, splitting monads are the same thing as adjoint retracts. In the latter reference it is moreover shown that direct sums and splitting monads suffice to admit all lax limits and all lax colimits. Combining all this, it thus follows that the principally generated modules, as a class of weights, describe precisely the lax (co)completion of $\Q$.

\section{Locally principally generated modules}\label{E}

It is well-known that idempotents in $\Sup$ split: for an idempotent $e\:L\to L$ we let
$$\xymatrix@=15mm{
L_e\ar@<1mm>[r]^{s_e} & L\ar@<1mm>[l]^{p_e} }$$
denote the obvious splitting with $L_e:=\{x\in L\mid e(x)=x\}$; of course any other splitting of $e$ is isomorphic to this one.
\par
Given a quantaloid $\Q$, let $\Q\si$ denote its split-idempotent completion; note that $\Q\si$ is small because $\Q$ is. Writing $j\:\Q\to\Q\si$ for the obvious inclusion, 
$$-\circ j\:\Mod(\Q\si)\to\Mod(\Q)$$
is an equivalence of quantaloids. We wish to describe the full subcategory of $\Mod(\Q)$ that is equivalent to $\Mod\pg(\Q\si)$ under the action of $-\circ j$. Thereto we shall first fix some notations.
\par
We shall write
$$(-)\si\:\Mod(\Q)\to\Mod(\Q\si)$$
for the inverse equivalence to $-\circ j$: it sends a $\Q$-module $\F$ to the $\Q\si$-module $\F\si$ defined (up to isomorphism) by:
\begin{itemize}
\item for an object $e\:A\to A$ of $\Q\si$, $\F\si(e):=\F(A)_{\F(e)}$,
\item for a morphism $f\:e\to e'$ in $\Q\si$, $\F\si(f):=\F(f)$.
\end{itemize}
\par
If $e\:A\to A$ is an idempotent in the quantaloid $\Q$, then the representable $\Sup$-natural transformation $\Q(-,e)\:\Q(-,A)\tto\Q(-,A)$ is idempotent too. All idempotents in $\Mod(\Q)$ split, so this one does too: we shall write
$$\xymatrix@=15mm{\F_{e}\ar@{=>}@<2mm>[r]^-{\sigma_e} & \Q(-,A)\ar@{=>}@<2mm>[l]^-{\pi_e}}$$
for the obvious splitting over $\F_e:=\Q\si(-,e)\circ j=\Q\si(j-,e)$, and we refer to such a $\Q$-module $\F_e$ as the {\em fixpoint $\Q$-module} for $e\:A\to A$.
\begin{proposition}\label{15}
Let $\F\:\Q\op\to\Sup$ be a $\Q$-module, and $\F\si$ the associated $\Q\si$-module.
\begin{enumerate}
\item\label{15.1} Given an idempotent $e\:A\to A$ in $\Q$, $a\in\F\si(e)$ is a principal element of the $\Q\si$-module $\F\si$ if and only if $a\in\F(A)$ satisfies $\F(e)(a)=a$ and the $\Q$-module morphism $\tau_a\circ\sigma_e\:\F_e\tto\F$ is a left adjoint.
\item\label{15.2} $\F\si$ is a principally generated $\Q\si$-module if and only if
\begin{equation}\label{16.0}
\mbox{for each $x\in\F$, }\tau_x=\bigvee\left\{\tau_a\circ[\tau_a,\tau_x]\bigmid 
\begin{array}{l}
\mbox{there exists an idempotent $e$} \\
\mbox{in $\Q$ such that $\F(e)(a)=a$ and} \\
\mbox{$\tau_a\circ\sigma_e\:\F_e\tto\F$ is a left adjoint}
\end{array}
\right\}.
\end{equation}
\end{enumerate}
\end{proposition}
\proof  
(\ref{15.1}) Because we put $\F\si(e)=\F(A)_{\F(e)}$, it is trivial that $a\in\F\si(e)$ if and only if $a\in\F(A)$ satisfying $\F(e)(a)=a$. By the $\Sup$-enriched Yoneda Lemma, we know that $a\in\F(A)$ corresponds uniquely to a $\Q$-module morphism, which we called $\tau_a\:\Q(-,A)\tto\F$; but also $a\in\F\si(e)$ corresponds uniquely to a $\Q\si$-module morphism: let us call it $\rho_a\:\Q\si(-,e)\tto\F\si$. By Definition \ref{12}, $a\in\F\si(e)$ is a principal element if $\rho_a$ is a left adjoint (in $\Mod(\Q\si)$, that is). Because of the equivalence of $\Mod(\Q)$ and $\Mod(\Q\si)$, expressed by $(-)\si$ and $-\circ j$, and because $\F_e=\Q\si(-,e)\circ j$, $\rho_a$ corresponds uniquely to a left adjoint $\Q$-module morphism $\zeta_a\:\F_e\tto\F$. But because $\F_e$ is the $\Q$-module over which the idempotent $\Q(-,e)\:\Q(-,A)\tto\Q(-,A)$ splits (and recall that we write the splitting with inclusion $\sigma_e$ and projection $\pi_e$) we necessarily have that $\zeta_a=\tau_a\circ\sigma_e$ (and $\tau_a=\zeta_a\circ\pi_e$). Thus it is indeed sufficient and necessary that $\tau_a\circ\sigma_e$ be a left adjoint $\Q$-module morphism.
\par
(\ref{15.2}) Again following Definition \ref{12}, and with the notations that we introduced in the first part of this proof, $\F\si$ is a principally generated $\Q\si$-module if
$$\mbox{for each $x\in\F\si$, }\rho_x=\bigvee\Big\{\rho_a\circ[\rho_a,\rho_x]\bigmid a\in(\F\si)\pr\Big\}.$$
This supremum of $\Q\si$-module morphisms can be written in terms of $\Q$-module morphisms, for similar reasons as in the first part of the proof:
$$\mbox{for each $x\in\F\si$, }\zeta_x=\bigvee\Big\{\zeta_a\circ[\zeta_a,\zeta_x]\bigmid a\in(\F\si)\pr\Big\}.$$
Using the notation
$$\boldalpha:=\left\{(a,e)\bigmid 
\begin{array}{l}
\mbox{$e\:A\to A$ an idempotent in $\Q$ and $a\in\F$ such } \\
\mbox{that $\F(e)(a)=a$ and $\tau_a\circ\sigma_e$ is left adjoint}
\end{array}\right\}$$
we can spell this out as:
\begin{equation}\label{16.2}
\begin{array}{l}
\mbox{for each idempotent $d\:X\to X$ in $\Q$} \\
\mbox{and each $x\in\F$ such that $\F(d)(x)=x$,}\\
\tau_x\circ \sigma_d=\bigvee\left\{(\tau_a\circ\sigma_e)\circ[(\tau_a\circ\sigma_e),(\tau_x\circ\sigma_d)]\bigmid 
(a,e)\in\boldalpha\right\}.
\end{array}
\end{equation}
\par
Assume now that \eqref{16.2} holds. Taking in particular $d=1_A$ (in which case $\sigma_d=\sigma_{1_A}$ is the identity transformation) it implies that for all $x\in\F$
\begin{eqnarray*}
\tau_x
 & = & \bigvee\left\{(\tau_a\circ\sigma_e)\circ[(\tau_a\circ\sigma_e),\tau_x]\bigmid (a,e)\in\boldalpha\right\} \\
 & = & \bigvee\left\{(\tau_a\circ\sigma_e)\circ\pi_e\circ[\tau_a,\tau_x]\bigmid (a,e)\in\boldalpha\right\} \\
 & = & \bigvee\left\{\tau_a\circ[\tau_a,\tau_x]\bigmid (a,e)\in\boldalpha\right\} \\
 & = & \bigvee\left\{\tau_a\circ[\tau_a,\tau_x]\bigmid a\in\boldbeta\right\}
\end{eqnarray*}
by a suitable application of Lemma \ref{17} stated below (and the proof of which is straightforward) to pass from the first to the second line, by the fact that $\tau_a\circ\sigma_e\circ\pi_e=\tau_a\circ\Q(-,e)=\tau_a$ (which is the equivalent of $\F(e)(a)=a$) to pass from the second to the third line, and where we introduced another auxiliary notation in the last line:
$$\boldbeta:=\Big\{a\bigmid\mbox{there exists an idempotent $e\:A\to A$ in $\Q$ such that }(a,e)\in\boldalpha\Big\}.$$
{\it A priori} an $a\in\F$ may be locally principal at several different idempotents, in which case $\boldalpha$ contains several ``copies'' of $a$ (one for each idempotent it is locally principal at) but $\boldbeta$ contains $a$ only once. But because an expression like $\tau_a\circ[\tau_a,\tau_x]$ does not contain any reference to the idempotents at which $a$ is locally principal, we can make the last step in the above series of equalities. Hence we derived the condition expressed in \eqref{16.0}.
\par
Conversely, assume the validity of \eqref{16.0}. Then, for every idempotent $d\:X\to X$ in $\Q$ and $x\in\F$ such that $\F(d)(x)=x$, we can compute in a similar fashion that
\begin{eqnarray*}
\tau_x\circ\sigma_d
 & = & \left(\bigvee\left\{\tau_a\circ[\tau_a,\tau_x]\bigmid a\in\boldbeta\right\}\right)\circ\sigma_d \\
 & = & \bigvee\left\{\tau_a\circ[\tau_a,\tau_x]\circ\sigma_d\bigmid a\in\boldbeta\right\} \\
 & = & \bigvee\left\{\tau_a\circ(\sigma_e\circ\pi_e)\circ[\tau_a,\tau_x]\circ\sigma_d\bigmid (a,e)\in\boldalpha\right\} \\
 & = & \bigvee\left\{(\tau_a\circ\sigma_e)\circ[(\tau_a\circ\sigma_e),(\tau_x\circ\sigma_d)]\bigmid (a,e)\in\boldalpha\right\}.
\end{eqnarray*}
Thus we obtain the condition expressed in \eqref{16.2}.
\endofproof
\begin{lemma}\label{17}
In any quantaloid, for a diagram like
$$\xymatrix@R=15mm@C=12mm{
X\ar@(ul,dl)_{p\circ i=1_X}\ar@<-1mm>@{>->}[d]_i & & Y\ar@<-1mm>@{>->}[d]_j\ar@(ur,dr)^{1_Y=q\circ j} \\
E\ar@(ul,dl)_{i\circ p=e^2=e}\ar@<-1mm>@{->>}[u]_p\ar[dr]_{a\circ e=a} & & F\ar@(ur,dr)^{f=f^2=j\circ q}\ar@<-1mm>@{->>}[u]_q\ar[dl]^{b=b\circ f} \\
 & A }$$
we have $q\circ[b,a]\circ i=[b\circ j,a\circ i]$.
\end{lemma}
\par
The result in Proposition \ref{15} suggests a new definition, to be compared with Definition \ref{12}:
\begin{definition}\label{20}
Let $\F\:\Q\op\to\Sup$ be a $\Q$-module. An $a\in\F$ is a {\em locally principal element (at an idempotent $e\:A\to A$ in $\Q$)} if (there is an idempotent $e\:A\to A$ in $\Q$ such that) $\F(e)(a)=a$ and $\tau_a\circ\sigma_e$ is a left adjoint in $\Mod(\Q)$. Writing $\F\lpr$ for the set of locally principal elements of $\F$, we say that $\F$ is {\em locally principally generated} if
$$\mbox{for each $x\in\F$, }\tau_x=\bigvee\Big\{\tau_a\circ [\tau_a,\tau_x]\bigmid a\in\F\lpr\Big\}.$$
\end{definition}
Thus, a locally principal element of $\F\:\Q\op\to\Sup$ {\it at an identity of $\Q$} is the same thing as, simply, a principal element of $\F$: {\it idempotents} in $\Q$ are viewed as {\it localities} (or {\it ``opens''}). It follows from the comparison of Definitions \ref{12} and \ref{20} that a principally generated $\Q$-module is necessarily also locally principally generated; but the converse is not true in general, as the next example shows.
\begin{example}\label{21.0}
A fixpoint $\Q$-module $\F_e\:\Q\op\to\Sup$ for an idempotent $e\:A\to A$ in $\Q$ trivially has $e\in\F_e(A)$ as locally principal element (at $e$, as a matter of fact); it follows straightforwardly that $\F_e$ is locally principally generated. However, $\F_e$ need not have any principal element and thus need not be principally generated; for a concrete example, let $\Q$ be the one-object suspension of the three-element chain $\{0\leq e\leq 1\}$.
\end{example}
\par
We have the following useful characterisation of locally principally generated $\Q$-modules.
\begin{proposition}\label{18.0}
Let $\F\:\Q\op\to\Sup$ be a $\Q$-module. The set of locally principal elements is
$$\F\lpr=\left\{\zeta(e)\bigmid 
\begin{array}{l}
\mbox{$e$ is an idempotent in $\Q$ and} \\
\mbox{$\zeta\:\F_e\tto\F$ is a left adjoint}
\end{array}
\right\},$$
and $\F$ is locally principally generated if and only if
\begin{equation}\label{16.1}
\id_{\F}=\bigvee\left\{\zeta\circ\zeta^*\bigmid 
\begin{array}{l}
\mbox{$e$ is an idempotent in $\Q$ and} \\
\mbox{$\zeta\:\F_e\tto\F$ is a left adjoint}
\end{array}
\right\}.
\end{equation}
\end{proposition}
\proof
This follows straightforwardly from Proposition \ref{11} and (the proof of) Proposition \ref{15}.
\endofproof
Much like Corollary \ref{11.3} does for principally generated $\Q$-modules, we can now give a characterisation in terms of universal constructions in $\Mod(\Q)$ of locally principally generated $\Q$-modules.
\begin{corollary}\label{18.1}
An $\F\in\Mod(\Q)$ is locally principally generated if and only if it is an adjoint retract of a direct sum of fixpoint $\Q$-modules.
\end{corollary}
\proof
Apply Lemma \ref{11.4} to the family of left adjoint $\Q$-module morphisms from all fixpoint $\Q$-modules to $\F$; compare with \eqref{16.1}.
\endofproof
\begin{remark}\label{21.1} 
A $\Q$-module $\F\:\Q\op\to\Sup$ is {\em projective} when the representable $\Sup$-functor $\Mod(\Q)(\F,-)\:\Mod(\Q)\to\Sup$ preserves epimorphisms. It is known (see e.g.\ [Stubbe, 2007a, Proposition 9.5]) that this is equivalent to $\F$ being {\em small-projective} (in the sense of [Kelly, 1982]: $\Mod(\Q)(\F,-)$ preserves all small weighted colimits), and equivalent to $\F$ being a retract of a direct sum of representable $\Q$-modules. It thus follows from Corollaries \ref{11.3} and \ref{18.1} that any (locally) principally generated $\Q$-module is necessarily projective in $\Mod(\Q)$. But the difference between ``projective'' and ``(locally) principally generated'' lies precisely in the strictly stronger requirement that, for the latter to hold, $\F$ needs to be an {\em adjoint} retract of a direct sum of representable $\Q$-modules (fixpoint $\Q$-modules).
\end{remark}
\par
Let $\Mod\lpg(\Q)$ denote the full subquantaloid of $\Mod(\Q)$ whose objects are the locally-principally generated $\Q$-modules.
\begin{theorem}\label{18}
For a small quantaloid $\Q$, the biequivalence $\Mod(\Q\si)\simeq\Mod(\Q)$ restricts to a biequivalence $\Mod\pg(\Q\si)\simeq\Mod\lpg(\Q)$.
\end{theorem}
The biequivalences in Corollary \ref{3.2} and Theorem \ref{13} apply to any small quantaloid, hence in particular to the split-idempotent completion $\Q\si$ of a small quantaloid $\Q$. The combination with the biequivalence in Theorem \ref{18} then shows that the following diagram commutes:
$$\xymatrix@R=8ex@C=4ex{
 & \Cocont(\Q\si)\ar[r]^-{\sim} & \Mod(\Q\si)\ar[r]^-{\sim} & \Mod(\Q) \\
\Cat\cc(\Q\si)\ar[r]^-{\sim}\ar[ru]^-{\P} & \Map(\Cocont\ta(\Q\si))\ar[r]^-{\sim}\ar@{^{(}->}[u] & \Map(\Mod\pg(\Q\si))\ar[r]^-{\sim}\ar@{^{(}->}[u] & \Map(\Mod\lpg(\Q))\ar@{^{(}->}[u]}$$
Up to the identification of $\Ord(\Q)$, the locally ordered category of ordered sheaves on a small quantaloid $\Q$, with the locally ordered category $\Cat\cc(\Q\si)$, the composition of the 2-functors
$$\Ord(\Q)=\Cat\cc(\Q\si)\arr{\P}\Cocont(\Q\si)\arr{\sim}\Mod(\Q\si)\arr{\sim}\Mod(\Q)$$
is precisely the left biadjoint 2-functor $F\:\Ord(\Q)\to\Mod(\Q)$ in \eqref{1} (see also [Stubbe, 2007b, 3.3]). Hence we proved that this 2-functor factors as a (composition of) biequivalence(s) followed by an inclusion.
\begin{theorem}\label{19}
The locally ordered category $\Ord(\Q)$ of ordered sheaves on a small quantaloid $\Q$ is biequivalent to $\Map(\Mod\lpg(\Q))$, the locally ordered category of locally principally generated $\Q$-modules and left adjoint $\Q$-module morphisms between them.
\end{theorem}

\section{Examples}\label{F}

In the three examples that follow we shall consider a one-object quantaloid $\Q$. In this situation we prefer to view $\Q$ as monoid $Q:=(\Q(*,*),\circ,1)$ in $\Sup$, and a $\Q$-module ${\cal M}\:\Q\op\to\Sup$ as the object $M:={\cal M}(*)$ of $\Sup$ together with an action\footnote{As mentioned earlier, $\Q$-modules are essentially ``the same thing'' as cocomplete $\Q$-categories; actions of $\Q$ correspond with tensors in $\Q$-categories. In previous work [Stubbe, 2006] we therefore denoted actions with ``$\tensor$'', the usual symbol for tensors in enriched categories. However to avoid any confusion with pure tensors in a tensor product of sup-lattices, we here adopt a ``$\circ$'' as notation.}
$$M\times Q\to M\:(x,q)\mapsto x\circ q:={\cal M}(q)(x).$$
With slight abuse of notation we shall write $[x,y]$ for the element of $Q$ that represents the lifting $[\tau_x,\tau_y]$ in $\Mod(Q)$; that is to say, $[x,y]=\bigvee\{q\in Q\mid x\circ q\leq y\}$ (compare with \eqref{6} and \eqref{11.0}). An element $a\in M$ is principal if and only if $[a,-]\:M\to Q$ is a $\Sup$-morphism that preserves the action of $Q$.
\par
Now suppose moreover that the underlying complete lattice of $Q$ is totally algebraic (in the classical sense, as recalled in the Introduction). Principality of an element $a\in M$ is then equivalent to the following two requirements:
\begin{enumerate}
\item $[a,x\circ c]\leq[a,x]\circ c$, for all $x\in M$ and all totally compact $c\in Q$,
\item $[a,\bigvee_{i\in I}x_i]\leq\bigvee_{i\in I}[a,x_i]$, for all $(x_i)_{i\in I}\in M$.
\end{enumerate}
This is equivalent to asking, for all $x,(x_i)_{i\in I}\in M$ and totally compact $c,d\in Q$,
\begin{enumerate}
\item if $a \circ d\leq x \circ c$ then there exists a totally compact $k\in Q$ such that $a \circ k\leq x$ and $d\leq k\circ c$,
\item if $a \circ d\leq\bigvee_{i\in I}x_i$ then there exists an $i\in I$ such that $a \circ d\leq x_i$.
\end{enumerate}
In particular, for a principal element $a\in M$ and a totally compact element $d\in\Q$, the element $a \circ d\in M$ is totally compact in (the underlying complete lattice of) $M$. Let $\Q\c$, resp.\ $M\c$, denote the partially orderd sets of totally compact elements of the underlying complete lattices of $\Q$, resp.\ $M$. For an element $x$ of a $\Q$-module $M$ we have
\begin{eqnarray*}
 & & x=\bigvee\Big\{a\circ[a,x]\bigmid a\in M\pr\Big\} \\
 & & \Longleftrightarrow \ x=\bigvee\Big\{a\circ q\bigmid a\in M\pr, q\in\Q,a\circ q\leq x\Big\} \\
 & & \Longleftrightarrow \ x=\bigvee\Big\{a\circ d\bigmid a\in M\pr, d\in\Q\c,a\circ d\leq x\Big\} \\
 & & \Longrightarrow \ x=\bigvee\Big\{b\bigmid b\in M\c,b\leq x\Big\}.
\end{eqnarray*}
Hence, every principally generated module on a totally algebraic quantale has a totally algebraic underlying sup-lattice.
\par
We shall now spell out three quite different applications.
\begin{example}[Complete lattices.] The quantaloid $\Sup$ can be identified with $\Mod(\2)$, where $\2$ stands for the two-element chain with its obvious quantale structure, $(\{0\leq 1\},\wedge,1)$. For a complete lattice/$\2$-module $S$, the conditions above say that:
\begin{enumerate}
\item $a\neq 0_S$,
\item $a$ is totally compact.
\end{enumerate}
Hence, the principal elements of $S$ are the non-zero totally compact elements, while $0_S$ is not principal but (the unique element that is) {\em locally} principal at $0$. Thus $S$ is a totally algebraic complete lattice if and only if it is a principally generated $\2$-module, if and only if it is a locally principally generated $\2$-module.
\end{example}
\begin{example}[Automata.] Let $(N,\cdot,1)$ be a monoid (in $\Set$), then the powerset of $N$ can be equiped with the pointwise multiplication
$$A\cdot B:=\{a\cdot b\mid a\in A, b\in B\}\mbox{ for $A,B\subseteq N$}$$
and thus $(2^N,\cdot,\{1\})$ is a quantale (it is the free quantale on $(N,\cdot,1)$, see [Rosenthal, 1990]). The totally compact elements of $2^N$ are the singletons and the empty subset. An element $a\in M$ of a $2^N$-module $M$ is principal if and only if, for all $x,(x_i)_{i\in I}\in M$ and $n,p\in N$,
\begin{enumerate}
\item [ia.] $a\circ n\neq 0_M$,
\item [ib.] if $a\circ n\leq x\circ p$ then there exists a $k\in N$ such that $a\circ k\leq x$ and $n=k\circ p$,
\item [ii.] if $a\circ n\leq\bigvee_{i\in I}x_i$ then there exists an $i\in I$ such that $a \circ n\leq x_i$.
\end{enumerate}
(We wrote $a\circ n$ instead of $a\circ\{n\}$ for notational convenience.) Note that (ia) + (ii) are equivalent to $a\circ n$ being principal in the underlying complete lattice of the module $M$.
\par
This example can be interesting for the theory of automata (or labelled transition systems): by Corollary \ref{13.1}, the principally generated $2^N$-modules can be identified with Cauchy complete $2^N$-enriched categories. It is well-known that categories enriched in a free quantale are precisely non-deterministic automata with $N$ as set of labels [Betti, 1980; Rosenthal, 1990].
\end{example}
\begin{example}[$\Rel(S,S)$-modules.] Let $\Rel$ denote the quantaloid of sets and relations, then surely for any set $S$, $Q_S:=\Rel(S,S)$ is a totally algebraic quantale: its totally compact elements are the empty set and the singletons $(s,t)\in S\times S$ (we omit the curly brackets for clarity). A $Q_S$-module $M$ is ``the same thing'' as the skeletal (i.e.\ having no non-identical isomorphic objects) cocomplete $Q_S$-category $\bbA_M$ (as explained in the beginning of Section \ref{D}). On the other hand, to give a $Q_S$-category $\bbA$, with object set $\bbA_0$, is equivalent to giving an order relation $\preccurlyeq$ on the set $\bbA_0\times S$; the correspondence is given by: 
$$(a,s)\preccurlyeq(b,t)\mbox{ if and only if }(s,t)\in\bbA(a,b).$$
Writing the equivalence relation on $\bbA_0\times S$ induced by the order $\preccurlyeq$ as $\approx$, it can be verified that $\bbA$ is skeletal and cocomplete if and only if $\preccurlyeq$ satisfies:
\begin{enumerate}
\item\label{x1} for all $s\in S$, $(a,s)\approx(b,s)$ implies $a=b$,
\item\label{x2} for all $(a_i)_{i\in I}\in\bbA_0$ there exists a $\bigvee_{i\in I}a_i\in\bbA_0$ such that 
$(\bigvee_{i\in I}a_i,s)\preccurlyeq(b,t)$ if and only if $(a_i,s)\preccurlyeq(b,t)$ for all $i\in I$,
\item\label{x3} for all $a\in\bbA_0$ and $s,t\in S$ there exists $b\in\bbA_0$ such that 
$(b,t)\approx(a,s)$ and moreover, whenever $u\neq t$, $(b,u)$ is a bottom element for the order $(\bbA_0\times S,\preccurlyeq)$.
\end{enumerate}
Conditions \eqref{x2} and \eqref{x3} are easily deduced from the equivalence of cocomplete $\Q$-categories to conically cocomplete and tensored $\Q$-categories [Kelly, 1982; Stubbe, 2006]: $\bigvee_{i\in I}a_i$ in \eqref{x2} is the conical colimit of $(a_i)_{i\in I}$ in $\bbA$ (and thus its order theoretical join in $(\bbA_0,\leq)$), while $b$ in \eqref{x3} is the tensor product $a\circ(s,t)$. Moreover, these conditions imply
\begin{itemize}
\item[ii'.] For all $(b_i)_{i\in I}\in\bbA_0$ there exists a $\bigwedge_{i\in I}b_i\in\bbA_0$ such that $(a,s)\preccurlyeq(\bigwedge_{i\in I}b_i,t)$ if and only if $(a,s) \preccurlyeq(b_i,t)$ for all $i\in I$.
\end{itemize}
Hence, a $Q_S$-module $M$ can be given in terms of an order $\preccurlyeq $ on $M\times S$, satifying conditions (\ref{x1}--\ref{x3}). An element $a\in M$ is then principal if and only if, for all $s,t\in S$, 
\begin{enumerate}
\item $(a,s)\not\approx\bot$, where $\bot$ denotes a bottom element for $(M\times S,\preccurlyeq)$,
\item if $(a,s)\preccurlyeq(\bigvee_{i\in I}x_i,t)$ then there exists an $i\in I$ such that $(a,s)\preccurlyeq(x_i,t)$.
\end{enumerate}
The order $\preccurlyeq $ on $M\times S$, corresponding to a $Q_S$-module $M$ has the characteristics of an entailment (especially \eqref{x2} and (ii') above): a couple $(a,s)\in M\times S$ can be thought of as an occurrence of an event $a\in M$ at a place $s\in S$. Quantales of the form $\Rel(S,S)$ arise in the context of relational representations of spatial quantales $\Q$, i.e.\ quantale homomorphisms $\rho\:\Q\to\End(2^S)=\Rel(S,S)$ [Mulvey and Resende, 2005].
\end{example}
\par
In the next section we shall dwell on the case where $\Q$ is the one-object suspension of a locale $X$. The formulation of ordered sheaves on $X$ by means of locally principally generated $X$-modules allows for a neat translation to ``skew local homeomorphisms'' into $X$.

\section{Skew local homeomorphisms}\label{SX}

\subsection*{Induced modules on a locale $X$}

In what follows, $\Loc$ denotes the (2-)category of locales. We follow the notational convention of [Johnstone, 1982, p.~40] for morphisms in $\Loc$: thus a locale morphism $f\:Y\to X$ is an adjoint pair
$$\xymatrix@=15mm{
Y\ar@{}[r]|{\bot}\ar@<-1mm>@/_2mm/[r]_{f_*} & X\ar@<-1mm>@/_2mm/[l]_{f^*}}$$
in the 2-category of partially ordered sets such that the left adjoint preserves finite infima. We do {\em not} follow the convention of [Johnstone, 1982; Mac~Lane and Moerdijk, 1992] when it comes to defining an order on the hom-sets in $\Loc$: for $f,g\:Y\biar X$ in $\Loc$ we define that $f\leq g$ if $f_*\leq g_*$. That is to say: we have that $\Loc\cong\Frm\coop$ as 2-categories (whereas the cited references have $\Loc\cong\Frm\op$)\footnote{The reason for our preference is in the first place notational convenience, especially in the 2-functors considered further on. However, there is maybe a deeper reason why this different ordering of locale morphisms is natural here: In the cited references locale morphisms are studied as inducing geometric morphisms between toposes of sheaves; the ordering of locale morphisms is chosen to correspond with the usual notion of natural transformation between geometric morphisms. We however shall study locale morphisms (or rather, morphisms in the slice category $\Loc\slice{X}$) as inducing order-preserving morphisms between the (ordered) sheaves themselves; and the ordering of the locale morphisms is chosen to correspond with the natural ordering of those morphisms between sheaves.}.
\par
Considering a locale $X$ as a monoid $(X,\wedge,\top_{\hspace{-0.6ex}X})$ in $\Sup$ it makes sense to write $\Mod(X)$ for the quantaloid of modules on the locale. Instead of writing these modules as contravariant $\Sup$-enriched presheaves on the one-object suspension of the locale, we rather consider them as objects of $\Sup$ on which $(X,\wedge,\top_{\hspace{-0.6ex}X})$ acts on the right: we write $(M,\circ)$ for a $\Sup$-object $M$ together with the action $(m,x)\mapsto m\circ x$. In the same vein, an $X$-module morphism $\alpha\:(M,\circ)\to(N,\circ)$ is a $\Sup$-morphism $\alpha\:M\to N$ which is equivariant for the respective actions.
\par
Given an $f\:Y\to X$ in $\Loc$, it is easily seen that putting
\begin{equation}\label{S5}
y\circ_f x:=y\wedge f^*(x)
\end{equation}
for $y\in Y$ and $x\in X$ results in an action of the monoid $(X,\wedge,\top_{\hspace{-0.6ex}X})$ on $Y$ in $\Sup$. In other words, from $f\:Y\to X$ in $\Loc$ we get an object $(Y,\circ_f)\in\Mod(X)$. Moreover, suppose that
\begin{equation}\label{S5.0}
\begin{array}{c}
\xymatrix{
Y\ar[rr]^h\ar[dr]_f & & Z\ar[dl]^g \\
 & X}
\end{array}
\end{equation}
is a commutative triangle in $\Loc$, then $h^*\:Z\to Y$ is a morphism in $\Sup$ satisfying $h^*(z\circ_f x)=h^*(z)\circ_g x$, for all $x\in X$, $z\in Z$. That is to say, $h^*\:(Z,\circ_g)\to(Y,\circ_f)$ is a morphism in $\Mod(X)$. All this adds up to an injective and faithful (but not full) 2-functor
\begin{equation}\label{S5.1}
(\Loc\slice{X})\coop\to\Mod(X).
\end{equation}
\par
We are now interested in left adjoint $X$-module morphisms:
\begin{definition}\label{S6}
A morphism $h\:f\to g$ in $\Loc\slice{X}$ as in \eqref{S5.0} is {\em skew open} if the corresponding order-preserving function $h^*\:Z\to Y$ has a left adjoint $h_!\:Y\to Z$ satisfying the ``balanced Frobenius identity\footnote{Putting $Z=X$ and $g=1_X$ this reduces to what is called the ``Frobenius identity'' in [Mac~Lane and Moerdijk, 1992, p.~500]; we call this generalisation  ``balanced'' because we get the (``unbalanced'') Frobenius identity by plugging in a terminal object.}'': for all $y\in Y$ and $x\in X$,
\begin{equation}\label{S7}
h_!(y\wedge f^*(x))=h_!(y)\wedge g^*(x).
\end{equation}
\end{definition}
\begin{example}\label{S6.0}
For an $h\:Y\to Z$ in $\Loc$ the following are equivalent:
\begin{enumerate}
\item $h\:Y\to Z$ is open in $\Loc$ (according to the ``usual'' definition of openness as in e.g.\ [Mac~Lane and Moerdijk, 1992, p.~500]),
\item for any $f\:Y\to X$ and $g\:Z\to X$ in $\Loc$ such that $g\circ h=f$, the morphism $h\:f\to g$ in $\Loc\slice{X}$ is skew open,
\item considering $h\:Y\to Z$ and $1_Z\:Z\to Z$ as objects in $\Loc\slice{Z}$, the (unique) morphism $h\:h\to 1_Z$ in $\Loc\slice{Z}$ is skew open.
\end{enumerate}
\end{example}
Clearly the identity morphisms in $\Loc\slice{X}$ are skew open, and the composition of skew open morphisms is again skew open; it thus makes sense to speak of the sub-2-category $(\Loc\slice{X})\o$ of $\Loc\slice{X}$ with the same objects but only its skew open morphisms. 
\par
Upon inspection it is easily seen that, for any two locale morphisms $f\:Y\to X$ and $g\:Z\to X$, there is an isomorphism of ordered sets
$$(\Loc\slice{X})\o(f,g)\cong\Map(\Mod(X))((Y,\circ_f),(Z,\circ_g))$$
given by sending a skew open morphism $h$ to the $X$-module morphism $h_!$ with right adjoint $h^*$. Sending skew open morphisms in $\Loc\slice{X}$ to their utmost left adjoints (i.e.\ $h\mapsto h_!$) thus gives rise to an injective and fully faithful 2-functor
\begin{equation}\label{S8}
(\Loc\slice{X})\o\to\Map(\Mod(X)).
\end{equation}
\par
In the codomain category of this functor we are now interested in the locally principally generated objects. In the next subsection we develop that notion further.

\subsection*{Locally principally generated $X$-modules}

Let $X$ be a locale. As is customary in locale theory, see e.g.\ [Mac~Lane and Moerdijk, 1992, p.~486], for any $u\in X$ we generically write $i\:\down u\mono X$ for the corresponding open sublocale of $X$, i.e.\ it is the open $\Loc$-morphism defined by $i_*(v):=(u\impl v)$, $i^*(x):=(x\wedge u)$ and $i_!(v):=v$. As noted before, it is therefore also skew open in $\Loc\slice{X}$ as (unique) morphism from $i\:\down u\mono X$ to the terminal object $1_X\:X\to X$,
\begin{equation}\label{S9.0}
\begin{array}{c}
\xymatrix{
\down u\ \ar@{>->}[dr]_i\ar@{>->}[rr]^i &  & X\ar[dl]^{1_X} \\
 & X }
\end{array}
\end{equation}
\par
All elements of the $\Sup$-monoid $(X,\wedge,\top_{\hspace{-0.6ex}X})$ are idempotent, thus each $u\in X$ gives rise to an idempotent representable $X$-module morphism on the (only) representable $X$-module $(X,\wedge)$. The image under the functor in \eqref{S8} of the $(\Loc\slice{X})\o$-morphism in \eqref{S9.0} is precisely the splitting of this idempotent:
\begin{equation}\label{S9}
\xymatrix@=15mm{
(\down u,\wedge)\ar@<1.5mm>[r]^{i_!} & (X,\wedge).\ar@<1.5mm>[l]^{i^*}\ar@(ul,ur)^{-\wedge u}}
\end{equation}
It is noteworthy that this is actually an {\em adjoint splitting}, since $i_!\dashv i^*$ in $\Mod(X)$, and that -- because $1_X$ is terminal in $\Loc\slice{X}$ and the functor in \eqref{S7} is fully faithful -- this is the {\em only} adjunction in $\Mod(X)$ between $(\down u,\wedge)$ and $(X,\wedge)$.
\par
Applying Definition \ref{20} to an $X$-module $(M,\circ)$ we get the following. An element $p\in M$ is {\em locally principal at $u\in X$} if and only if $p\circ u=p$ and the composite $X$-module morphism
$$(\down u,\wedge)\arr{i_!}(X,\wedge)\arr[9mm]{p\circ-}(M,\circ)$$
admits a right adjoint in $\Mod(X)$. Let $(M,\circ)\lpr$ denote the set of elements of $M$ which are locally principal at some $u\in X$. Then $(M,\circ)$ is {\em locally principally generated} if and only if, for each $m\in M$, 
$$m=\bigvee\{p\circ[p,m]\mid p\in(M,\circ)\lpr\},$$
where $[p,m]:=\bigvee\{u\in X\mid p\circ u\leq m\}$.
\par
We shall recast the latter definition in a more pleasant form.
\begin{proposition}\label{S11}
Let $(M,\circ)$ be an $X$-module.
\begin{enumerate}
\item\label{S11.a} If $p\in M$ is locally principal at $u\in X$, then for any $m\in M$, $p\circ[p,m]$ is locally principal at $u\wedge[p,m]$.
\item\label{S11.b} For any $m\in M$, $\{p\circ[p,m]\mid p\in(M,\circ)\lpr\}=\down m \cap (M,\circ)\lpr.$
\item\label{S11.c} $(M,\circ)$ is locally principally generated if and only if
\begin{equation}\label{S12}
\mbox{for all $m\in M$, }m=\bigvee(\down m \cap (M,\circ)\lpr).
\end{equation}
\end{enumerate}
\end{proposition}
\proof
\eqref{S11.a} For shorthand we introduce $q:=p\circ[p,m]$ and $v:=u\wedge[p,m]$. Then it is easily verified that
\begin{eqnarray*}
q\circ v
 & = & (p\circ[p,m])\circ(u\wedge[p,m]) \\
 & = & p\circ([p,m]\wedge u\wedge[p,m]) \\
 & = & (p\circ u)\circ[p,m] \\
 & = & p\circ[p,m] \\
 & = & q.
\end{eqnarray*}
Moreover the diagram
$$\xymatrix{
(\down v,\wedge)\ar[r]^{i_!}\ar[d] & (X,\wedge)\ar[r]^{q\circ-} & (M,\circ)\ar@{=}[d] \\
(\down u,\wedge)\ar[r]_{i_!} & (X,\wedge)\ar[r]_{p\circ-} & (M,\circ)}$$
in $\Mod(X)$, where the left downward arrow is the obvious inclusion of $\down v$ into $\down u$, commutes: for $w\leq v$ we can compute that
$$q\circ w=(p\circ[p,m])\circ w=p\circ([p,m]\wedge w)=p\circ w.$$
But $(\down v,\wedge)\to(\down u,\wedge)\:w\mapsto w$ is a left adjoint in $\Mod(X)$, hence the top composite morphism is a left adjoint whenever the bottom composite morphism is.
\par
\eqref{S11.b} Because $p\circ-\dashv[p,-]$ as order-preserving maps between $X$ and $M$, it is trivial that $p\circ[p,m]\leq m$ and ($p\leq m\impl p\circ[p,m]=p$), for any $p,m\in M$. We have just shown that if $p$ is locally principal then so is $p\circ[p,m]$. Hence the equality of these sets.
\par
\eqref{S11.c} Is now immediate.
\endofproof
\par
We can also translate to an $X$-module $(M,\circ)$ the condition in Proposition \ref{18.0} that expresses that it is locally principally generated if and only if
\begin{equation}\label{S12.0}
\id_{(M,\circ)}=\bigvee\{\zeta\circ\zeta^*\mid u\in X, \zeta\in\Map(\Mod(X))((\down u,\wedge),(X,\wedge))\},
\end{equation}
where we write $\zeta^*$ for the right adjoint to $\zeta$. This fact allows us to prove the following remarkable property:
\begin{proposition}\label{S12.1}
Let $(M,\circ)$ be a locally principally generated $X$-module.
\begin{enumerate}
\item\label{S12.a} For $m,n\in M$, $m=n$ if and only if for every $u\in X$ and every left adjoint $X$-module morphism $\zeta\:(\down u,\wedge)\to(M,\circ)$ we have $\zeta^*(m)=\zeta^*(n)$.
\item\label{S12.b} For every $m,n\in M$ and $x\in X$, $(m\wedge n)\circ x=m\wedge(n\circ x)$.
\item\label{S12.c} $M$ is a locale and $f^*\:X\to M\:x\mapsto \top_{\hspace{-.6ex}M}\circ x$
is the inverse image of a locale morphism $f\:M\to X$ for which $(M,\circ_f)=(M,\circ)$.
\end{enumerate}
\end{proposition}
\proof
\eqref{S12.a} One direction is trivial; for the other one expands $m=\id_{(M,\circ)}(m)$ and $n=\id_{(M,\circ)}(n)$ by means of the formula given above.
\par
\eqref{S12.b} Let $u\in X$ and $\zeta\in\Map(\Mod(X))((\down u,\wedge),(X,\wedge))$ with right adjoint $\zeta^*$. Then $\zeta^*((m\wedge n)\circ x)=\zeta^*(m\wedge n)\wedge x=(\zeta^*(m)\wedge\zeta^*(n))\wedge x$ because $\zeta^*$ is a module morphism (and thus turns the ``$-\circ x$'' into a ``$-\wedge x$'') and because it is a right adjoint (and thus preserves infima). But ``for the same reasons'' we also have that $\zeta^*(m\wedge(n\circ x))=\zeta^*(m)\wedge\zeta^*(n\circ x)=\zeta^*(m)\wedge(\zeta^*(n)\wedge x))$. Thus $\zeta^*((m\wedge n)\circ x)=\zeta^*(m\wedge(n\circ x))$ for all $u$ and all $\zeta$, and we conclude by the above that $(m\wedge n)\circ x=m\wedge(n\circ x)$.
\par
\eqref{S12.c} Let $m,(m_i)_{i\in I}$ be elements of $M$. Let $u\in X$ and $\zeta\:(\down u,\wedge)\to(X,\wedge)$ a left adjoint in $\Mod(X)$ with right adjoint $\zeta^*$. Using that $\zeta^*$ is both a left and a right adjoint in $\Ord$ one computes that
$$\zeta^*(m\wedge\bigvee_i m_i)=\zeta^*(m)\wedge\bigvee_i\zeta^*(m_i)$$
but also that 
$$\zeta^*(\bigvee_i(m\wedge m_i))=\bigvee_i(\zeta^*(m)\wedge\zeta^*(m_i)).$$
In both right hand sides we now find elements of the locale $\down u$, where $\wedge$ distributes over $\bigvee$, and hence these expressions are equal. This holds for all $u$ and all $\zeta$, so by the the first statement we obtain $m\wedge\bigvee_i m_i=\bigvee_i(m\wedge m_i)$, which means that $M$ is a locale. Finally, the function $f^*\:X\to M$ is certainly a $\Sup$-morphism: because the action of $(X,\wedge,\top_{\hspace{-0.6ex}X})$ on $M$ preserves suprema ``in both variables''. But moreover, for $x,y\in X$, we may compute -- using the formula in \eqref{S12.b} with $m=\top_M\circ x$ and $n=\top_M$ to pass from the second line to the third -- that
\begin{eqnarray*}
f^*(x\wedge y)
 & = & \top_{\hspace{-0.6ex}M}\circ(x\wedge y) \\
 & = & (\top_{\hspace{-0.6ex}M}\circ x)\circ y \\
 & = & (\top_{\hspace{-0.6ex}M}\circ x)\wedge(\top_{\hspace{-0.6ex}M}\circ y) \\
 & = & f^*(x)\wedge f^*(y).
\end{eqnarray*}
Thus $f^*$ is indeed the inverse image part of a locale morphism $f\:Y\to X$. Putting $n=\top_{\hspace{-0.6ex}M}$ in the formula in \eqref{S12.b} it follows that moreover
$$m\wedge f^*(x)=m\wedge(\top_{\hspace{-0.6ex}M}\circ x)=m\circ x,$$
that is to say, $(M,\circ_f)=(M,\circ)$ as claimed.
\endofproof
\par
We now go on to define the notion of ``skew local homeomorphism''.

\subsection*{Skew local homeomorphisms}

Let $f\:Y\to X$ be in $\Loc$ and $u\in X$; we keep the notation $i\:\down u\mono X$ for the corresponding open sublocale of $X$. Recall from [Mac~Lane and Moerdijk, 1992, p.~524] that the elements of the set
$$S_f(u):=\Loc\slice{X}(i,f)$$
are the {\em sections of $f$ at $u$}. This defines a sheaf $S_f\:X\op\to\Set$, and this construction extends to a functor $\Loc\slice{X}\to\Sh(X)$ whose restriction to local homeomorphisms is an equivalence of categories.
\par
A particular feature of local homeomorphisms is that, whenever $f=g\circ h$ in $\Loc$, if $f$ and $g$ are local homeomorphisms then so is $h$; recall also that a local homeomorphism is always open in $\Loc$ (see {\it loc.\ cit.}). Thus, if $f\:Y\to X$ is a local homeomorphism then every $s\in S_f(u)$ is an {\em open section} in the sense that $s\:\down u\to Y$ is an open locale morphism. With this in mind the following is a natural generalisation.
\begin{definition}\label{S13}
For $f\:Y\to X$ in $\Loc$ and $i\:u\mono X$, we put
$$S_f\o(u):=(\Loc\slice{X})\o(i,f)$$
and call its elements the {\em skew open sections of $f$ at $u$}.
\end{definition}
\begin{example}\label{S14}
Every open section $s\:\down u\to Y$ of a locale map $f\:Y\to X$ is necessarily skew open too; but the converse need not hold. However, if $f\:Y\to X$ is a local homeomorphism then $S_f(u)=S_f\o(u)$ for all $u\in X$.
\end{example}
\par
A morphism $f\:Y\to X$ in $\Loc$ is a local homeomorphism if and only if $Y$ can be covered by its open sections [Johnstone, 2002, vol.~2, p.~503], i.e.\
$$\top_Y=\bigvee\left\{s_!(u)\bigmid 
\mbox{$u\in X$, $s\in S_f(u)$ and $s$ is open in $\Loc$}\right\}.$$ 
In this case, every $y\in Y$ can be covered by open sections of $f$, by taking the restrictions of the open sections of $f$ to $y$. This motivates our main definition in this section:
\begin{definition}\label{S17}
A morphism $f\:Y\to X$ in $\Loc$ is a {\em skew local homeomorphism} if 
$$1_Y=\bigvee\{s_!\circ s^*\mid u\in X, s\in S_f\o(u)\}.$$
\end{definition}
For the record we immediately add:
\begin{example}\label{S17.1}
Every local homeomorphism is a skew local homeomorphism. A skew local homeomorphism is a local homeomorphism if and only if its (skew open) sections are all open.
\end{example}
Skew local homeomorphisms can be characterised in different ways:
\begin{proposition}\label{S15}
Let $f\:Y\to X$ be in $\Loc$. 
\begin{enumerate}
\item There is a bijection between skew open sections of $f$ at $u\in X$ and locally principal elements of the $X$-module $(Y,\circ_f)$ at the idempotent $u\in (X,\wedge,\top_{\hspace{-0.6ex}X})$; if $s\in S_f\o(u)$ then $s_!(u)\in Y$ is the corresponding locally principal element.
\item\label{S16.0} The following statements are equivalent:
\begin{enumerate}
\item $(Y,\circ_f)$ is a locally principally generated $X$-module,
\item for all $y\in Y$, $y=\bigvee(\down y\cap \{s_!(u)\mid u\in X, s\in S_f\o(u)\})$,
\item $f$ is a skew local homeomorphism.
\end{enumerate}
\end{enumerate}
\end{proposition}
\proof
(i) By the fully faithful 2-functor in \eqref{S8} we know, for each $i\:\down u\mono X$, that
$$(\Loc\slice{X})\o(i,f)\cong\Map(\Mod(X))((\down u,\wedge),(Y,\circ_f));$$
the left hand side is precisely $S_f\o(u)$, and the bijection is given from left to right by sending an $s\in S_f\o(u)$ to $s_!$. As in Proposition \ref{18.0}, the right hand side is in bijection with the set of elements of $(Y,\circ_f)$ which are locally principal at $u$, by sending $s_!$ to $s_!(u)$.
\par
(ii) Immediate from \eqref{S11.c} in Proposition \ref{S11}, and \eqref{S12.0}.
\endofproof
\par
Let $(\Loc\slice{X})\slh\o$ denote the full subcategory of $(\Loc\slice{X})\o$ whose objects are the skew local homeomorphisms. It follows from the above results that the fully faithful 2-functor in \eqref{S8} (co)restricts to a fully faithful 2-functor
\begin{equation}\label{S18}
(\Loc\slice{X})\slh\o\to\Map(\Mod\lpg(X)).
\end{equation}
This 2-functor is easily seen to be injective on objects; but due to Proposition \ref{S12.1} it is surjective too: for every locally principally generated $X$-module $(M,\circ)$ the locale morphism $f\:M\to X$ with inverse image $f^*(x)=\top_{\hspace{-0.6ex}M}\circ x$, which satisfies $(M,\circ)=(M,\circ_f)$, is a skew local homeomorphism. The consequence of our work is then the following result.
\begin{theorem}\label{S20}
For any locale $X$, the 2-functor in \eqref{S18} is an isomorphism of locally ordered categories:
$$(\Loc\slice{X})\o\slh\cong\Map(\Mod\lpg(X));$$
both of these are thus equivalent to $\Ord(X)\simeq\Ord(\Sh(X))$, the ordered sheaves on $X$ viewed as enriched categorical structures, resp.\ the internal orders in the topos $\Sh(X)$.
\end{theorem}
\par
We have seen in Example \ref{S17.1} that any local homeomorphism is necessarily a skew local homeomorphism; and we have seen in Example \ref{S6.0} that any open locale morphism is necessarily skew open too. It follows that $\LH\slice{X}$ is a full subcategory of $(\Loc\slice{X})\o\slh$. We just proved the latter to be isomorphic to $\Map(\Mod\lpg(X))$, thus it makes sense to determine those locally principally generated $X$-modules which, under this isomorphism, correspond to local homeomorphisms.
\begin{definition}\label{S22}
A locally principally generated $X$-module $(M,\circ)$ is an {\em \'etale $X$-module} when every left adjoint $X$-module morphism $\zeta\:(\down u,\wedge)\to(M,\circ)$ satisfies, for all $v\in\down u$ and $m\in M$,
$$\zeta(v\wedge\zeta^*(m))=\zeta(v)\wedge m.$$
\end{definition}
It is straightforward from Example \ref{S17.1} and Proposition \ref{S15} that a skew local homeomorphism $f\:Y\to X$ is a local homeomorphism if and only if $(Y,\circ_f)$ is an \'etale $X$-module. Letting $\Mod\et(X)$ stand for the full sub-2-category of $\Mod\lpg(X)$ consisting of \'etale $X$-modules, we can conclude with the following summary.
\begin{theorem}\label{S23}
For any locale $X$ there is a commuting square
$$\xymatrix@=15mm{
(\Loc\slice{X})\o\slh\ar@{=}[r] & \Map(\Mod\lpg(X)) \\
\LH\slice{X}\ar@{=}[r]\ar@{^{(}->}[u] & \Map(\Mod\et(X))\ar@{^{(}->}[u]}$$
in which the equalities denote isomorphisms of (locally ordered) categories and the upward arrows are full embeddings. The categories in the bottom row are equivalent to $\Sh(X)$, the locally ordered categories in the top row are equivalent to $\Ord(\Sh(X))$, and the inclusions view ``sets as discrete (or symmetric) orders''.
\end{theorem}

\section{Addendum: biadjunction, biequivalence}\label{AA}

It was shown in [Stubbe, 2007a, Section 8] that, for any $\Q$-category $\bbC$, the totally compact objects of the presheaf category $\P\bbC$ form precisely the Cauchy-completion of $\bbC$: $(\P\bbC)\c=\bbC\cc$. That is to say, a contravariant presheaf $\phi\:*_A\to\bbC$ is totally compact in $\P\bbC$ if and only if it has a right adjoint in $\Dist(\Q)$ (``$\phi$ is Cauchy''). Representable contravariant presheaves certainly are Cauchy, thus the Yoneda embedding $Y_{\bbC}\:\bbC\to\P\bbC$, which sends an object $c\in\bbC$ to the representable $\bbC(-,c)\:*_{tc}\dist\bbC$, corestricts to $(\P\bbC)\c$:
\begin{equation}\label{A5.1}
\overline{Y}_{\bbC}\:\bbC\to(\P\bbC)\c\:c\mapsto\bbC(-,c).
\end{equation}
The following is a mere triviality.
\begin{lemma}\label{A5.1.1}
For any $\Q$-category $\bbC$, the functor $\overline{Y}_{\bbC}\:\bbC\to(\P\bbC)\c$ is fully faithful, and it is an equivalence if and only if $\bbC$ is Cauchy-complete.
\end{lemma}
\par
On the other hand we can compute, for any cocomplete $\Q$-category $\bbA$ and any $\phi\in\P(\bbA\c)$, the $\phi$-weighted colimit of the inclusion $i_{\bbA}\:\bbA\c\to\bbA$. This defines a functor
\begin{equation}\label{A5.2}
R_{\bbA}\:\P(\bbA\c)\to\bbA\:\phi\mapsto\colim(\phi,i_{\bbA})
\end{equation}
about which we record some auxiliary results.
\begin{lemma}\label{A6}
For any cocomplete $\Q$-category $\bbA$, the functor $R_{\bbA}\:\P(\bbA\c)\to\bbA$ is cocontinuous and admits a cocontinuous right adjoint. Moreover, $R_{\bbA}$ is always fully faithful, and it is an equivalence if and only if $\bbA$ is totally algebraic.
\end{lemma}
\proof
We claim that $R_{\bbA}$ is left adjoint to $H:\bbA\to\P(\bbA\c):a\mapsto\bbA(i_{\bbA}-,a)$. Indeed, for $\phi\in\P(\bbA\c)$ and $x\in\bbA$,
\begin{eqnarray*}
\bbA(R_{\bbA}\phi,x)
 & = & \bbA(\colim(\phi,i_{\bbA}),x) \\
 & = & \Big[\phi,\bbA(i_{\bbA}-,x)\Big] \\
 & = & \P(\bbA\c)(\phi,Hx).
\end{eqnarray*}
\par
Next we prove that $H$ itself is cocontinuous; it suffices to show that it preserves suprema of contravariant presheaves: for any $\phi\in\P\bbA$, $H(\sup_{\bbA}(\phi))=\colim(\phi,H)$. Note first that, by Proposition \ref{4},
\begin{eqnarray*}
H(\sup_{\bbA}(\phi)) 
 & = & \bbA(i_{\bbA}-,\sup_{\bbA}(\phi)) \\
 & = & \phi(i_{\bbA}-) \\
 & = & \bbA(i_{\bbA}-,-)\tensor\phi.
\end{eqnarray*}
But then also
\begin{eqnarray*}
\P(\bbA\c)\Big(H(\sup_{\bbA}(\phi)), -\Big)
 & = & \Big[\bbA(i_{\bbA}-,-)\tensor\phi,-\Big] \\
 & = & \Big[\phi,[\bbA(i_{\bbA}-,-),-]\Big] \\
 & = & \Big[\phi,\P(\bbA\c)(H-,-)\Big]
\end{eqnarray*}
which is the universal property that we had to check.
\par
To see that $R_{\bbA}$ is fully faithful it is (necessary and) sufficient to show that the unit of the adjunction $R_{\bbA}\dashv H$ is an isomorphism (cf.\ [Stubbe, 2007, 2.3] for example). Thus, for $c\in\bbA\c$ and $\phi\in\P(\bbA\c)$ we compute that
\begin{eqnarray*}
((H\circ R_{\bbA})(\phi))(c) 
 & = & \bbA(i_{\bbA}c,\colim(\phi,i_{\bbA})) \\
 & = & H_c(\colim(\phi,i_{\bbA})) \ \ \ \mbox{(with notations as in Proposition \ref{4})} \\
 & = & \colim(\phi,H_c\circ i_{\bbA}) \ \ \ \mbox{(because of Proposition \ref{4})} \\
 & = & \colim(\bbA(-,i_{\bbA}-)\tensor\phi,H_c).
\end{eqnarray*}
By an argument in the proof of Proposition \ref{4}, we know that a presheaf-weighted colimit of $H_c$ is the value of that weight in $c$; here this allows us to equate 
$$\colim(\bbA(-,i_{\bbA}-)\tensor\phi,H_c)=\bbA(c,i_{\bbA}-)\tensor\phi=\bbA\c(c,-)\tensor\phi=\phi(c),$$
taking into account that $c\in\bbA\c$ and $i_{\bbA}$ is a full embedding. This indeed proves that $H\circ R_{\bbA}=1_{\P(\bbA\c)}$.
\par
Finally, knowing that $R_{\bbA}$ is always fully faithful, it is an equivalence if and only if also the counit of the adjunction  $R_{\bbA}\dashv H$ is an isomorphism. Spelled out this means that, for every $a\in\bbA$, $a\cong\colim(\bbA(i_{\bbA}-,a),i_{\bbA})$, which precisely says that $1_{\bbA}$ is the (pointwise) left Kan extension of $i_{\bbA}$ along itself.
\endofproof
\begin{theorem}\label{A7}
There is a biadjunction
$$\Cat(\Q)\inlineadj{\P}{\ (-)\c}\Map(\Cocont(\Q))$$
where the involved 2-functors are defined as:
$$\begin{array}{rcl}
\P\Big(\bbC\arr{F}\bbD\Big) & := & \P\bbC\arr[25mm]{\bbD(-,F-)\tensor-}\P\bbD,\\
\Big(\bbA\arr{F}\bbB\Big)\c & := & \bbA\c\arr{F}\bbB\c
\end{array}$$
and counit and unit are given by the functors in \eqref{A5.1} and \eqref{A5.2}.
\end{theorem}
\proof
We shall prove that, for any cocomplete $\Q$-category $\bbA$,
$$\Big(\bbA\c,R_{\bbA}\:\P(\bbA\c)\to\bbA\Big)$$
is a biuniversal right reflection along the 2-functor $\P\:\Cat(\Q)\to\Map(\Cocont(\Q))$. The latter is indeed a 2-functor: for an $F\:\bbC\to\bbD$ in $\Cat(\Q)$ we have adjoints
$$\P\bbC
\xymatrix@=40mm{
\ar@{}@<8.5mm>[r]|{\bot}\ar@{}@<-8.5mm>[r]|{\bot}\ar@<3mm>@/^10mm/[r]^{\bbD(-,F-)\tensor-}\ar@<-3mm>@/_10mm/[r]_{[\bbD(F-,-),-]}
 & \ar[l]_{[\bbD(-,F-),-]}^{=\bbD(F-,-)\tensor-}}
\P\bbD\mbox{\ \ \ in }\Cat(\Q)$$
(with all compositions and liftings computed in $\Dist(\Q)$) so $\P$ lands in the 2-category $\Map(\Cocont(\Q))$; and 2-functoriality is obvious. Moreover, Lemma \ref{A6} provides the information that $R_{\bbA}$ is a morphism of $\Map(\Cocont(\Q))$. So all we need to show, is that $R_{\bbA}$ has the required 2-universal property, i.e.\ the order-preserving function
\begin{equation}\label{A8}
\Cat(\Q)(\bbC,\bbA\c)\to\Map(\Cocont(\Q))(\P\bbC,\bbA)\:F\mapsto R_{\bbA}\circ\P F
\end{equation}
is an equivalence of ordered sets. We shall prove first that it is essentially surjective, and then that it is order-reflecting.
\par
{\it Essential surjectivity.} Suppose given a left adjoint $G\:\P\bbC\to\bbA$ in $\Cocont(\Q)$, or equivalently, suppose given adjoints
$$\P\bbC\xymatrix@=15mm{
\ar@{}@<4mm>[r]|{\bot}\ar@{}@<-4mm>[r]|{\bot}\ar@<1.5mm>@/^5mm/[r]^{G}\ar@<-1.5mm>@/_5mm/[r]_{K}
 & \ar[l]|{H}}
\bbA\mbox{\ \ \ in }\Cat(\Q).$$
For $\psi\in\P\bbA$ we can compute, with straightforward arguments involving liftings and compositions in the quantaloid $\Dist(\Q)$, that $\psi(G\circ Y_{\bbC}-)$ is the $\psi$-weighted colimit of $H$:
\begin{eqnarray*}
\Big[\psi,\P\bbC(H-,-)\Big]
 & = & \Big[\psi,[\P\bbC(Y_{\bbC}-,H-),-]\Big] \\
 & = & \Big[\psi,[\bbA(G\circ Y_{\bbC}-,-),-]\Big] \ \ \ \mbox{(because $G\dashv H$)} \\
 & = & \Big[\bbA(G\circ Y_{\bbC}-,-)\tensor\psi,-\Big] \\
 & = & \Big[[\bbA(-,G\circ Y_{\bbC}-),\psi],-\Big] \\
 & = & \Big[\psi(G\circ Y_{\bbC}-),-\Big] \ \ \ \mbox{(by Yoneda Lemma for $\Q$-cats)} \\
 & = & \P\bbC(\psi(G\circ Y_{\bbC}-),-).
\end{eqnarray*}
But for any $x\in\bbC$, we can also compute that
\begin{eqnarray*}
\bbA(GY_{\bbC}(x),\sup_{\bbA}(\psi))
 & = & \P\bbC(Y_{\bbC}(x),H(\sup_{\bbA}(\psi))) \ \ \ \mbox{(by adjunction $G\dashv H$)} \\
 & = & (H\circ\sup_{\bbA}(\psi))(x) \ \ \ \mbox{(by Yoneda Lemma for $\Q$-cats)} \\
 & = & \colim(\psi,H)(x) \ \ \ \mbox{(because $H$ is cocontinuous)}.
\end{eqnarray*}
Putting these together we have, for any $\psi\in\P\bbA$ and $x\in\bbC$, that
$$\psi(GY_{\bbC}(x))=\bbA(GY_{\bbC}(x),\sup_{\bbA}(\psi))$$
which, according to Proposition \ref{4}, means that for any $x\in\bbC$ the object $GY_{\bbC}(x)$ of $\bbA$ is totally compact. In other words, the given $G\:\P\bbC\to\bbA$ factors as
$$\xymatrix@=15mm{
\P\bbC\ar[rr]^G & & \bbA \\
\bbC\ar[u]^{Y_{\bbC}}\ar@{.>}[r]_{\overline{G}} & \bbA\c\ar[ur]_i}$$
where $\overline{G}(x):=G(Y_{\bbC}(x))$. It is a matter of calculations, using cocontinuity of $G$ amongst other things, to see that $G=R_{\bbA}\circ\P(\overline{G})$: for $\phi\in\P\bbC$, 
\begin{eqnarray*}
(R_{\bbA}\circ \P(\overline{G}))(\phi)
 & = & \colim(\P(\overline{G})(\phi),i_{\bbA}) \\
 & = & \colim(\bbA\c(-,\overline{G}-)\tensor\phi,i_{\bbA}) \\
 & = & \colim(\bbA(-,i_{\bbA}-)\tensor\bbA\c(-,\overline{G}-)\tensor\phi,1_{\bbA}) \\
 & = & \colim(\bbA(-,i_{\bbA}\circ\overline{G}-)\tensor\phi,1_{\bbA}) \\
 & = & \colim(\bbA(-,GY_{\bbC}-)\tensor\phi,1_{\bbA}) \\
 & = & \colim(\phi,G\circ Y_{\bbC}) \\
 & = & G\circ\colim(\phi,Y_{\bbC}) \\
 & = & G(\phi).
\end{eqnarray*}
Thus we proved that the function in \eqref{A8} is essentially surjective.
\par
{\it Order-reflection.} Remark first that for any $F\:\bbC\to\bbA\c$ in $\Cat(\Q)$, the outer diagram
$$\xymatrix@=15mm{
\P\bbC\ar[r]\ar[r]^{\P(F)} & \P(\bbA\c)\ar[r]^{R_{\bbA}} & \bbA \\
\bbC\ar[u]^{Y_{\bbC}}\ar@{.>}[r]_{F} & \bbA\c\ar[u]^{Y_{\bbA\c}}\ar[ur]_i}$$
commutes: both the left hand square and the right hand triangle are easily checked by computation. Now suppose that some $F,G\:\bbC\biar\bbA\c$ in $\Cat(\Q)$ are such that $R_{\bbA}\circ \P(F)\leq R_{\bbA}\circ\P(G)$ in $\Map(\Cocont(\Q))$. Then we can deduce from the above that
$$i_{\bbA}\circ F = R_{\bbA}\circ\P(F)\circ Y_{\bbC}\leq R_{\bbA}\circ\P(G)\circ Y_{\bbC}=i_{\bbA}\circ G.$$
But because $i_{\bbA}\:\bbA\c\to\bbA$ is a full embedding, it follows that necessarily $F\leq G$ from the start. This proves that the function in \eqref{A8} is also order-reflecting.
\par
It is now a matter of routine computations to verify that the right biadjoint to the 2-functor $\P\:\Cat(\Q)\to\Map(\Cocont(\Q))$ is indeed given by ``restricting to totally compacts'':
$$(-)\c\:\Map(\Cocont(\Q))\to\Cat(\Q)\:\Big(F\:\bbA\to\bbB\Big)\mapsto\Big(F\:\bbA\c\to\bbB\c\Big),$$
and that the unit of the biadjunction is indeed given by those corestrictions of Yoneda embeddings as in \eqref{A5.1}.
\endofproof
The (co)restriction of biadjoint 2-functors to those objects for which the (co)unit is an equivalence, is a biequivalence of 2-categories. In the case of interest above, we recover via Lemmas \ref{A5.1.1} and \ref{A6} the biadjunction in Corollary \ref{3.2}.

\end{document}